\documentclass[icol, sn-mathphys-ay]{sn-jnl}

\usepackage{graphicx}%
\usepackage{multirow}%
\usepackage{amsmath,amssymb,amsfonts}%
\usepackage{amsthm}%
\usepackage{mathrsfs}%
\usepackage[title]{appendix}%
\usepackage{xcolor}%
\usepackage{textcomp}%
\usepackage{manyfoot}%
\usepackage{booktabs}%
\usepackage{algorithm}%
\usepackage{algorithmicx}%
\usepackage{algpseudocode}%
\usepackage{listings}%
\usepackage{hyperref}       
\usepackage{url}            
\usepackage{nicefrac}       
\usepackage{microtype}      
\usepackage{stmaryrd}
\usepackage{geometry}

\newcommand{\R}{\mathbb{R}}
\newcommand{\N}{\mathbb{N}}
\newcommand{\E}{\mathbb{E}}
\newcommand{\V}{\mathbb{V}}

\newcommand{\tr}{\text{Tr}}

\newcommand{\Cov}{\text{Cov}}

\usepackage{physics}

\DeclareMathOperator*{\minimize}{minimize}

\theoremstyle{thmstyleone}%
\newtheorem{thm}{Theorem}
\newtheorem{ppt}{Proposition}%

\newtheorem{lem}{Lemma}

\theoremstyle{thmstyletwo}%
\newtheorem{rmk}{Remark}

\theoremstyle{thmstylethree}%
\newtheorem{dft}{Definition}
\newtheorem{as}{Assumption}
\newtheorem{notat}{Notation}

\begin{document}

\title[Article Title]{Analysis of a multi-target linear shrinkage covariance estimator}

\author*[1,2]{\fnm{Benoît} \sur{Oriol}}\email{benoit.oriol@dauphine.eu}

\affil*[1]{\orgdiv{CEREMADE}, \orgname{Université Dauphine-PSL}, \orgaddress{\street{16, Boulevard de Lannes}, \city{Paris}, \postcode{75116}, \country{France}}}

\affil[2]{\orgdiv{Core ML Lab}, \orgname{Société Générale Corporate and Investment Banking}, \orgaddress{\street{4, Boulevard des Bouvets}, \city{Puteaux}, \postcode{92800}, \country{France}}}

\abstract{Multi-target linear shrinkage is an extension of the standard single-target linear shrinkage for covariance estimation. We combine several constant matrices - the targets - with the sample covariance matrix. We derive the oracle and a \textit{bona fide} multi-target linear shrinkage estimator with exact and empirical mean. In both settings, we proved its convergence towards the oracle under Kolmogorov asymptotics. Finally, we show empirically that it outperforms other standard estimators in various situations.}


\maketitle

\section*{Acknowledgments} 
We want to give special thanks to Alexandre Miot (Core ML  Lab, Société Générale Corporate and Investment Banking) and Gabriel Turinici (CEREMADE, Université Dauphine-PSL) for their advices all along the work.

\section{Introduction and related work}

The covariance matrix plays a major role in numerous machine learning algorithms and statistics. Just to cite a few, array signal processing \cite{Abrahamsson2007, Guerci1999}, generalized method of moments estimators \cite{Hansen1982}, financial portfolio management \cite{Markowitz1952, Ledoit2003} and forecasting \cite{Tsay2009}, functional genomics \cite{Schafer2005}. However, those algorithms are designed to use the population covariance matrix, which is often unaccessible. The sample covariance matrix severely fails in many applications: for instance, for Markowitz portfolio management it does often not beat a naive uniform diversification among the assets \cite{Demiguel2009}.

In the context of Kolmogorov asymptotics, where the ratio of the dimension $p_n$ and the number of samples $n$ tends to a finite positive constant $\frac{p_n}{n} \rightarrow c > 0$, the sample covariance eigenvalue spectrum is spread: high eigenvalues tends at being too high, and low ones, too low \cite{Girko1985, Girko1990, Girko1995}.

Influenced by the work of Stein on Gaussian mean estimation in 1956 \cite{Stein1956}, we focus on a specific type of regularization: shrinkage estimation. The motivation is to design a covariance estimator that behaves well in high dimension, with a lower variance than the sample covariance at the cost of an increased bias. This variance-bias trade-off aims at reducing the overall mean-squared error (MSE) of the covariance estimator with respect to the population covariance $\Sigma$. Linear shrinkage estimators linearly combine the sample covariance matrix with a finite set of predefined constant matrices - the targets. 

More formally, we denote by $S$ the sample covariance, and we consider $N \in \N$ targets, denoted $(T^{(i)})_{i=1}^N$, which are constant matrices of the same size as $S$. Usually, when there is only one target, $N=1$, we chose $T^{(1)} = I$ the identity matrix \cite{Ledoit2004}. With $N=2$ targets, a standard example is the following. Consider a diagonal block population covariance $
	\Sigma = \begin{pmatrix}
		\Sigma_1 & (0) \\
		(0) & \Sigma_2
	\end{pmatrix}.$ The method does \textbf{not} require $\Sigma$ to be block-diagonal, this example uses it only for intuition and explainatory reasons. Even if $\Sigma_1$ and $\Sigma_2$ are unknown, we can leverage our knowledge of the block structure, and chose the following targets: 
\begin{equation*}
\begin{aligned}
	T^{(1)} = \begin{pmatrix}
		I_1 & (0) \\
		(0) & (0)
	\end{pmatrix} \text{ and } 
	T^{(2)} = \begin{pmatrix}
		(0) & (0) \\
		(0) & I_2
	\end{pmatrix}.
\end{aligned}
\end{equation*}
where $I_1$ and $I_2$ are the identity matrix of same size as, respectively, $\Sigma_1$ and $\Sigma_2$. A linear shrinkage covariance estimator $\hat \Sigma$ is, by definition:
\begin{equation*}
\begin{aligned}
	\hat \Sigma = c_0 S + \sum_{i=1}^N c_i T^{(i)} \text{, where } (c_0,...,c_N) \in \R^N.
\end{aligned}
\end{equation*}
It is a linear combination of the sample covariance and the targets. The challenge is to estimate $(c_0,...,c_N)$ only with the observations in order to minimize the Frobenius norm $\lVert \hat \Sigma - \Sigma \rVert_F^2$, where $\Sigma$ is the population covariance.

Many single-target shrinkage estimators (STSE) were proposed, just to name a few: Ledoit-Wolf linear shrinkage \cite{Ledoit2004}, its extension for Gaussian distributions named Oracle Approximating Shrinkage (OAS) estimator \cite{Chen2010}, linear shrinkage with factor models \cite{Ikeda2016a}, linear shrinkage for elliptical distributions with unknown mean and known radius distribution \cite{Ashurbekova2021}. 

However, the performance of those algorithms highly depends on the choice of the target and its similarity to the population covariance $\Sigma$. Thus, there is no general ideal choice of target fitting all covariances $\Sigma$. The use of several targets instead of one can lead to superior classification performance \cite{Halbe2013, Lancewicki2014}. The methods to compute multi-target linear shrinkage estimators are less numerous. With only two targets, an analytic estimator is given in \cite{Halbe2013}. The general case for any number of targets an estimator was proposed for Gaussian distributions in \cite{Gray2018}. In the distribution-free case, an implicit analytic estimator was proposed in \cite{Lancewicki2014}, retrieved from solving a strictly convex problem. To the best of our knowledge, there is no theoretical analysis of those algorithms comparable to the quadratic convergences we have in single-target shrinkage. 

We address the lack of explicit estimators and theoretical results for multi-target shrinkage covariance estimation. This paper proposes an explicit analytic distribution-free multi-target shrinkage estimator (MTSE). It provides an asymptotic theoretical analysis, proving its quadratic convergence to the oracle estimator in the scope of Kolmogorov asymptotics. Experimentally, our estimator outperforms single-target estimators in the sense of MSE in various situations.

The paper is organized as follows. In section 2, we introduce the notation, definitions and hypotheses necessary to define and address the MTSE problem. In section 3, we expose the problem, the explicit oracle and proposed estimator, and show its quadratic convergence. Section 4 presents the numerical simulations. Finally, section 5 exposes the principal limitations of the approach.

For easier reading, technical proofs of the theoretical results are derived in the appendix, we will only give the idea of proof for intuition in the main corpus. The code and data for the experiment is available on our \href{https://github.com/cvcer/MTSE}{GitHub}.

\section{Notation, definitions and hypotheses}
Let us introduce the following notation.

\begin{notat}[Observations]
We consider a sequence of observation matrices $(X_n)_{n \in \N^*} $ with $X_n \in \R^{p_n \times n}$ of $n$ iid observations on a system of $p_n$ dimensions. Decomposing the covariance matrix, we denote $\Sigma_n = \Gamma_n \Lambda_n \Gamma_n^t$, where $\Lambda_n$ is a diagonal matrix and $\Gamma_n$ a rotation matrix. The diagonal elements of $\Lambda_n$ are the eigenvalues $\lambda_1^n,...,\lambda_{p_n}^n$, and the columns of $\Gamma_n$ are the eigenvectors $\gamma_1^n,...,\gamma_{p_n}^n$. $Y_n = \Gamma_n^t X_n$ is a $p_n \times n$ matrix of $n$ iid observations of $p_n$ uncorrelated random variables $(y_1^n,...,y_n^n)$.
\end{notat}

\begin{notat}[Norm]
Let $A_n$ and $B_n$ be two $p_n \times p_n$ matrices. We consider the Frobenius norm: $ \lVert A_n \rVert_n = \sqrt{\tr(A_nA_n^T)/p_n}$, and the associated inner product: $\langle A_n, B_n \rangle_n = \tr(A_nB_n^T)/p_n$. Dividing by the dimension is not standard, it is done to fix the norm of the identity as $1$ regardless of the dimension.
\end{notat}

\begin{notat}[Convergence]
Let $(E_n)_n$ be a sequence of euclidean spaces with associated norm $\lVert \cdot \rVert_n$. The quadratic convergence of a random variable $Z_n \in E_n$, i.e. $\E[\lVert Z_n \rVert_n^2] \rightarrow 0$, is denoted as $Z_n \underset{q.m}{\longrightarrow} 0$.
\end{notat}

We define the sample covariance in the setting of known and unknown mean, and the set of targets.
\begin{dft}[Empirical covariance]
As we are studying both situations of known or unknown mean, we will use a different sample covariance $S_n$ depending on the situation. If we know the mean of the distribution, for an observation matrix $X_n$ of size $p_n \times n$, we define the empirical covariance as:
\begin{equation}
\begin{aligned}
	S_n = (X_n - \E[X]) (X_n - \E[X])^T/n.
\end{aligned}
\end{equation}
Otherwise, if we do not know the mean of the distribution, we will use the unbiased sample covariance:
\begin{equation}
\begin{aligned}
	&S_n = \tilde X_n \tilde X_n^T/(n-1),
\end{aligned}
\end{equation}
with $(\tilde X_n)_{ik} = (X_n)_{ik} - \frac{1}{n}\sum_{k'=1}^n (X_n)_{ik'}$.
\end{dft}

\begin{dft}[Set of targets]
For each $n$, we define the finite sequence of $N_n \in \N^*$ targets: $\left(T_n^{(k)} \right)_{k=1}^{N_n}$, where each $T_n^{(k)}$ is a symmetric matrix of size $(p_n,p_n)$, not necessarily positive. We assume that for each $n$, $\left(T_n^{(k)} \right)_{k=1}^{N_n}$ form a linearly independent family, orthogonal and normalized.
\end{dft}

The ``orthogonal and normalized'' condition is not really a restriction in practice. Indeed, we can build any linearly independent family $T_n^{(k)}$ of symmetric matrices. Then, we can use the Gram-Schmidt algorithm \cite{Gram1883} on this family to make it ortho-normalized, without modifying $\text{Span}\left(\left(T_n^{(k)} \right)_{k=1}^{N_n}\right)$, so without changing the space of multi-target estimators of the form $\hat \Sigma = c_0 S + \sum_{i=1}^N c_i T^{(i)} \text{, where } (c_0,...,c_N) \in \R^N$.

So, we consider without loss of generality that the family $\left(T_n^{(k)} \right)_{k=1}^{N_n}$ is orthogonal and normalized.


We describe now several assumptions, extensions to our framework of the assumptions used in the linear shrinkage of Ledoit and Wolf \cite{Ledoit2004}, which will be used in the following. The first one is the core hypothesis in Kolmogorov asymptotics: the dimension can grow along with the number of samples.
\begin{as}[Dimension]
There exists a constant $K_1$ independent of $n$ such that $p_n/n \leq K_1$.
\end{as}

Bounded $8$-th moments are necessary for the proofs, even if experimentally we will see that bounded $4$-th moments seem sufficient for the convergence to hold. For the notation, $y_{i1}^n$ designates the $i$-th dimension of the first sample in $Y_n$. As all samples are i.i.d., it could be any other arbitrary choice of sample in $\llbracket 1,n \rrbracket$ of course. Similarly, $X_{1}^n$ is the first sample, of size $p_n$, from $X_n$.
\begin{as}[Moments]
There exists a constant $K_2$ independent of $n$ such that $\frac{1}{p_n} \sum_{i=1}^{p_n} \E\left[(y_{i1}^n)^8\right] \leq K_2$ and for all $k \in \llbracket 1, N_n \rrbracket$, $\E\left[\left\langle X_{1}^n (X_{1}^n)^T, T_n^{(k)} \right\rangle_n^4\right] \leq K_2$.
\end{as}

The number of target $N_n$ can depend of $n$, while staying negligible compared to $n$. Intuitively, each target comes with an additional coefficient to estimate, and the estimation error decreases the performance of the \textit{bona fide} estimator.
\begin{as}[Number of targets]
The number of targets $N_n$ verifies that: $N_n = o(n)$.
\end{as}

Already present in \cite{Ledoit2004}, this technical assumption, not really restrictive in practice as argued in the same article, is still necessary here.
\begin{as}[Technical covariance assumption]
With $Q_n$ the set of all the quadruples that are made of four distinct integers between $1$ and $p_n$, and for all $i \in \llbracket 1,p_n \rrbracket$, $\tilde y_{i1}^n = y_{i1}^n - \E[y_{i1}^n]$, we have:
\[\lim_{n \rightarrow \infty} \frac{p_n^2}{n} \times \frac{\sum_{(i,j,k,l) \in Q_n} (\Cov[\tilde y_{i1}^n \tilde y_{j1}^n, \tilde y_{k1}^n \tilde y_{l1}^n])^2}{|Q_n|} = 0.\]
\end{as}

Finally, for the MTSE problem to have a.s. a unique solution, we need the following characterization. In practice, if $\lVert S_n \rVert_n^2 - \sum_{k} \langle S_n, T_n^{(k)} \rangle_n^2 = 0$, \textit{i.e.} if $S_n \in \text{Span}\left(\left(T_n^{(k)} \right)_{k=1}^{N_n}\right)$, simply return $S_n$.
\begin{as}[Well-posed problem]\label{wellposed}
For each $n \in \N^*$, using the associated probability $\mathbb{P}_n$ of the distribution indexed by $n$,
\[\mathbb{P}_n\left[\lVert S_n \rVert_n^2 - \sum_{k} \langle S_n, T_n^{(k)} \rangle_n^2 = 0\right ] = 0.\]
\end{as}

\section{Oracle and proposed estimator}
Given a sample covariance $S_n$, we are looking for a minimizer of the 2-norm in the scope of linear combinations of $S_n$ and $T_n^{(k)}$, \textit{i.e.} the orthogonal projection of the population covariance onto the subspace spanned by $\{S_{n},T_{n}^{(1)},...,T_{n}^{(N_{n})}\}$.
\begin{equation}\label{min_or}
\begin{aligned}
	\minimize_{\substack{c \in \R^{N_n+1}}} \left \lVert c_0 S_n + \sum_{k=1}^{N_n} c_{k} T_n^{(k)} - \Sigma_n \right \rVert_n^2.
\end{aligned}
\end{equation} 

The idea is to solve this minimization problem, construct an \textbf{oracle} covariance estimator with it, and then estimate its components and propose a \textbf{bona fide} covariance estimator. Firstly, let us see the \textbf{oracle} covariance estimator $\Sigma_n^*$ that emerges.
\begin{ppt}[Oracle estimator]
The oracle coefficients solving the minimization problem \ref{min_or} are, for $i \in \llbracket 1, N_n \rrbracket$, with $\langle \cdot \rangle := \langle \cdot \rangle_n$:
\begin{equation}
\begin{aligned}
&c_0^* = \frac{1}{d_n^2}\left(\langle S_n, \Sigma_n \rangle - \sum_{k=1}^{N_n} \langle S_n, T_n^{(k)} \rangle\langle \Sigma_n, T_n^{(k)} \rangle \right),\\
&c_i^* =  \langle \Sigma_n,  T_n^{(i)} \rangle - c_0^* \langle S_n,  T_n^{(i)} \rangle, \\
&d_n^2 = \lVert S_n \rVert_n^2 - \sum_{k=1}^{N_n} \langle S_n, T_n^{(k)} \rangle_n^2.
\end{aligned}
\end{equation}
Thus, we define the \textbf{oracle} covariance estimator $
\Sigma_n^{*} = \mathcal{P}_{\mathcal{S}_n^+} \left(c_0^* S_n + \sum_{k=1}^{N_n} c_{k}^* T_n^{(k)}\right)$,
where $\mathcal{P}_{\mathcal{S}_n^+}$ is the orthogonal projection on the cone of positive semi-definite matrices.
\end{ppt}

Remark that even the formulas tend to look rough, the oracle estimator is none else but the \textbf{orthogonal projection} of the population covariance $\Sigma_n$ onto the subspace spanned by $\{S_{n},T_{n}^{(1)},...,T_{n}^{(N_{n})}\}$, followed by projection onto the cone of positive semi-definite matrices $\mathcal{S}_n^+$.

This orthogonal projection $\mathcal{P}_{\mathcal{S}_n^+}$ does not strictly follow the definition of linear shrinkage. However, as $\Sigma_n$ is symmetric positive, the projection step decreases the distance $\lVert \Sigma_n^* - \Sigma_n \rVert_n$. This step is not necessary for the convergence results to hold, but it is a free improvement we recommend.

Nevertheless, the oracle estimator still depends on $\Sigma_n$, and cannot be computed as it is. We propose a \textit{bona fide} estimator, and show that asymptotically it converges to the oracle one. The strategy is the following, both for known and unknown mean:
\begin{itemize}
	\item For $i \in \llbracket 1, N_n \rrbracket$, in $c_i^*$, $\langle \Sigma_n,  T_n^{(i)} \rangle_n$ is estimated by $\langle S_n,  T_n^{(i)} \rangle_n$,
	\item $\langle S_n, \Sigma_n \rangle_n$ is approximated by $\lVert S_n \rVert_n^2 - \underbrace{\E\left[\lVert S_n - \Sigma_n \rVert_n^2 \right]}_{\V(S_n)}$. Note that in expectation both objects are equal. $\V(S_n)$ is estimated by $\hat \V(S_n)$ defined below,
	\item Similarly, for each $k \in \llbracket 1, N_n \rrbracket$, $\langle S_n, T_n^{(k)} \rangle_n\langle \Sigma_n, T_n^{(k)} \rangle_n$ is approximated by $\langle S_n, T_n^{(k)} \rangle_n^2 - \underbrace{\E\left[\langle S_n - \Sigma_n, T_n^{(k)} \rangle_n^2\right]}_{\V\left(\langle S_n , T_n^{(k)} \rangle_n \right)}$. $\V\left(\langle S_n , T_n^{(k)} \rangle_n \right)$ is estimated by $\hat \V\left(\langle S_n , T_n^{(k)} \rangle_n \right)$, defined below.
\end{itemize}
\begin{dft}[Known mean - Variance estimators]
We suppose the mean to be \textbf{known}. We define:
\begin{equation*}
\begin{aligned}
	&\hat{\V}(S_n) = \frac{1}{n(n-1)} \sum_{k=1}^n \left \lVert X_{\cdot k}^n \left(X_{\cdot k}^n\right)^T - S_n \right \rVert_n^2.
\end{aligned}
\end{equation*}
Let $i \in \llbracket 1, N_n \rrbracket$. We define:
\begin{equation*}
\begin{aligned}
	&\hat{\V}(\langle S_n, T_n^{(i)} \rangle_n) = \sum_{k=1}^n \frac{\left \langle X_{\cdot k}^n \left(X_{\cdot k}^n\right)^T - S_n, T_n^{(i)} \right \rangle_n^2}{n(n-1)}.
\end{aligned}
\end{equation*}
\end{dft}
For the \textbf{unknown mean} setting, closed-form formulas are derived for $\hat{\V}(S_n)$ and $\hat{\V}(\langle S_n, T_n^{(i)} \rangle_n)$, and are detailed in the Appendix.
\begin{thm}[Variance estimators]
	Denoting $\V(S_n) = \E[\lVert S_n - \Sigma_n \rVert_n^2]$ and $\V\left(\langle S_n , T_n^{(i)} \rangle_n \right) = \E\left[\langle S_n - \Sigma_n, T_n^{(i)} \rangle_n^2\right] $, under assumptions 1 to 4, both where the mean is known or unknown for all $i \in \llbracket 1, N_n \rrbracket$, $\hat{\V}(S_n)$ and $\hat{\V}(\langle S_n, T_n^{(i)} \rangle_n)$ are unbiased estimators of $\V(S_n)$ and $\V(\langle S_n, T_n^{(i)} \rangle_n)$ respectively. Moreover, $\E\left[\left(\hat{\V}(S_n) - \V(S_n)\right)^2\right] \longrightarrow 0$ and $ \E\left[\left(\sum_{i=1}^{N_n} \hat{\V}(\langle S_n, T_n^{(i)} \rangle_n) - \V(\langle S_n, T_n^{(i)} \rangle_n\right)^2\right] \longrightarrow 0$.
\end{thm}

The idea of proof of this theorem resides in carefully developing, counting null terms and using the assumptions 3 and 4 through Cauchy-Schwarz and Jensen inequalities.

The previous results gave us \textit{bona fide} estimators - that behave well asyptotically - for each quantity necessary to compute $\Sigma_n^*$. Then, we can define our \textit{bona fide} covariance estimator $S_n^*$. As we carefully chose estimators that behave well asymptotically, we prove that $S_n^*$ converges to $\Sigma_n^*$ in quadratic mean.

\begin{dft}[Bona fide estimator]
The bona fide estimators $c$ for the minimization problem \ref{min_or} are defined, for $i \in \llbracket 1, N_n \rrbracket$, as:
\begin{equation}
\begin{aligned}
&c_0 = \min\left((\hat{c}_0)_+, 1\right), c_{i} =  (1 - c_0)\langle S_n,  T_n^{(i)} \rangle_n,
\end{aligned}
\end{equation}
\begin{equation*}
\begin{aligned}
\text{with }\hat{c}_{0} = &\frac{1}{d_n^2}\Bigg(\lVert S_n \rVert_n^2 - \hat{\V}(S_n) 
- \sum_{k=1}^{N_n} \left(\langle S_n, T_n^{(k)} \rangle_n^2 - \hat{\V} \left(\langle S_n, T_n^{(k)} \rangle_n\right) \right) \Bigg),\\
d_n^2 = &\lVert S_n \rVert_n^2 - \sum_{k=1}^{N_n} \langle S_n, T_n^{(k)} \rangle_n^2,
\end{aligned}
\end{equation*}
The \textit{bona fide} MTSE is defined as: 
\begin{equation*}
\begin{aligned}
S_n^* =\mathcal{P}_{\mathcal{S}_n^+} \left( c_0 S_n + (1-c_0) \sum_{i=1}^{N_n} c_{i}T_n^{(i)} \right).
\end{aligned}
\end{equation*}
\end{dft}
The implementation in Python is available on our \href{https://github.com/cvcer/MTSE}{GitHub}. The truncation of $c_0$ between $0$ and $1$ theoretically does not change the asymptotic behavior of the estimator. Geometrically, it ensures that $S_n^*$ is a convex combination between $S_n$ and the matrix $\sum_{i=1}^{N_n} c_{i}T_n^{(i)}$. Empirically, it improves the robustness of the estimator, particularly in small sample size.


\begin{thm}[Loss convergence]
	Under assumptions 1 to 5, $S_n^*$ quadratically converges to the optimal estimator $\Sigma_n^{*}$, i.e. $S_n^* - \Sigma_n^{*} \underset{q.m}{\longrightarrow} 0$. Consequently, $S_n^*$ has asymptotically the same loss as $\Sigma_n^{*}$, i.e. $\E\left[\left|\left\lVert S_n^* - \Sigma_n \right \rVert_n^2 - \left\lVert \Sigma_n^{*} - \Sigma_n \right \rVert_n^2 \right|\right] \rightarrow 0$.
\end{thm}

The idea of the proof is to prove the convergence on the estimators without the projection step onto the cone of positive matrices, and use the Lipschitz property of projection to conclude. The initial convergence relies on the Lemma A.1 in \cite{Ledoit2004}.

\section{Experimental results}
This section illustrates the practical consequences of the choice of targets on the estimation, the robustness of the estimator to heavy tails, and an application to real non-stationary time series. The experiments were computed on a Lenovo laptop, Intel i5, 32 GB Ram. For reproducibility, the code and data are available on our \href{https://github.com/cvcer/MTSE}{GitHub}, and the settings of the experiments are detailed in the Appendix. The performance is measured with the Percentage Relative Improvement in Average Loss (PRIAL):
\begin{equation}
\begin{aligned}
\text{PRIAL}(\hat{S}_n) = \frac{\E[\lVert S_n - \Sigma_n \rVert_n^2] - \E[\lVert \hat{S}_n - \Sigma_n \rVert_n^2]}{\E[\lVert S_n - \Sigma_n \rVert_n^2]}. 
\end{aligned}
\end{equation}
It is to be noted that the complexity of the implementation can be reduced to $O(N_np_n^3 + N_nnp_n)$. The details are in the Appendix.
\subsection{Target set impact: Usefulness of a target}
We fix the dimension at $p = 50$, and $n = 25$, and we study the impact of the number of targets $N_n$. The population covariance $\Sigma$ is designed to clearly exhibit a choice of useful targets, and compare the role of targets improve significantly the oracle PRIAL-, to those which give less. $\Sigma$ is made as follows - a formal definition is given in the Appendix -, where blocks $B_i$ are drawn in a Wishart distribution:
\begin{equation}\label{eq_block}
\begin{aligned}
	\Sigma = \begin{pmatrix}
		B_1 &  & (0) \\
		 & \ddots &   \\
		(0) & & B_{10}
	\end{pmatrix}.
\end{aligned}
\end{equation}

We consider two sets of targets: one set of ``aligned'' targets $\{T_1^{(i)}\}$, and one set of ``misaligned'' ones $\{T_2^{(i)}\}$. The first one is built as following: the first target $T_1^{(1)} = I$, to match Ledoit-Wolf STSE \cite{Ledoit2004}, and for $i \in \llbracket 2,10 \rrbracket$, the target $T_1^{(i)}$ is block-diagonal; where $\Sigma$ has the block $B_i$, $T_1^{(i)}$ has an Identity block:
\begin{equation} \label{eq_target}
\begin{aligned}
	T_1^{(i)} = \begin{pmatrix}
		(0) &  \ldots & (0) \\
		\vdots& (\text{block }i:) I_5   & \vdots \\
		(0)  & \ldots &  (0)
	\end{pmatrix}
\end{aligned}
\end{equation}

For the second set of targets $\{T_2^{(i)}\}$, we break the perfect alignment between $T_1^{(i)}$ and the block $B_i$. For that, we permute the rows and columns: column $j$ in $\{T_1^{(i)}\}$ becomes column $j+2$ in $\{T_2^{(i)}\}$ and so on.

\begin{figure}
\centering
\includegraphics[width=0.9\linewidth]{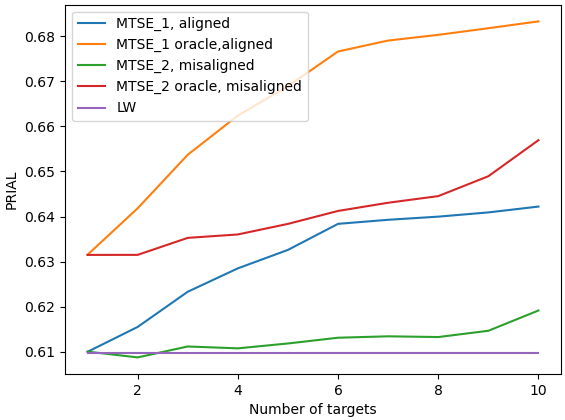}
\caption{PRIAL of the MTSE, and the oracle, with target sets 1 and 2, and the Ledoit-Wolf STSE as reference.}
\label{fig_offset}
\end{figure}

We consider the PRIAL of the MTSE (MTSE\_1 and MTSE\_2 using respectively $\{T_1^{(i)}\}$ and $\{T_2^{(i)}\}$), the oracle ones, and the Ledoit-Wolf STSE (LW), in function of the number $k$ of targets used, between $k=1$, which is equivalent to the STSE, and $k=10$ where the sets of targets are used completely. The results, with $50000$ draws and Student distribution ($\nu=9$ d.o.f.), are shown in Figure \ref{fig_offset}. The maximum standard deviation is $4 \times 10^{-4}$, error bars were removed.

\begin{figure}
	\centering
\includegraphics[width=0.9\linewidth]{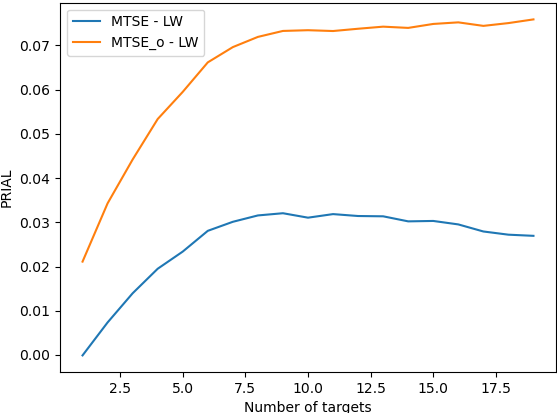}
\caption{PRIAL difference between the oracle MTSE and Ledoit Wolf estimator (MTSE\_o - LW), and bona fide MTSE (MTSE - LW).}
\label{fig:1}
\end{figure}


It confirms the intuition about a useful target: the ``aligned'' set $(T_1)^{(i)}$ presents a significant better PRIAL than the second misaligned one $(T_2)^{(i)}$. This comprehension of what makes a target useful is of utter importance to design the target set in a practical use case. The second conclusion is that even if the misaligned targets are not optimal, they still present a better PRIAL than the STSE using only one target.

\subsection{Target set impact: Adding useless targets to a good target set}
We keep the previous setting for $p$, $n$, $\Sigma$ and the samples $X$ but we change the set of targets.

The first target is the identity matrix $I$, to match with Ledoit-Wolf STSE. For  $i \in \llbracket 2,10 \rrbracket$, $T^{(i)}$ is designed according to Equation \ref{eq_target}. For target $i \geq 11$, the target is drawn from a Wishart distribution $\mathcal{W}_p(I_p,p)$. The latter are supposed to be ``useless'' regarding the structure of $\Sigma$.

We ran the experiment for $5000$ times, the maximum standard deviation is $7 \times 10^{-4}$. We plot the \textbf{differences} of PRIAL between the MTSE and the Ledoit-Wolf STSE, the same for the MTSE oracle. The result can be seen in figure \ref{fig:1}. 

Two distinct behaviors emerge from this experiment. The first $10$ useful targets results in a significant improvement of the PRIAL. For the random targets, they can be tagged as useless as they do not improve significantly the oracle PRIAL. The interesting conclusion is that the estimator PRIAL is only slightly negatively affected but stays above the STSE PRIAL. Consequently, \textbf{adding good targets gives a high reward, and adding useless targets comes at a low price}.

\begin{figure}
\centering
\includegraphics[width=0.9\linewidth]{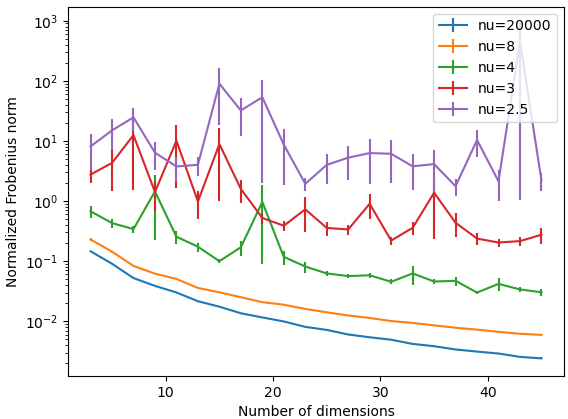}
\caption{Convergence of $\lVert S_n^* - \Sigma_n^{*}\rVert_n^2$ for different values of $\nu$.}
\label{fig:nu}
\end{figure}

\subsection{Heavy-tails robustness}
In the case of Ledoit-Wolf STSE, Ledoit and Wolf empirically argue that Assumption 2 about finite eighth moments is in practice not necessary, only finite fourth moments would suffice \cite{Ledoit2004}. Thus, we study the convergence of $S_n^*$ towards $\Sigma_n^{*}$ in a setting of heavy tails. For that, we consider multivariate t-Student distributions with $\nu > 2$ degrees of freedom, for different values of $\nu \in \{20000, 8, 4, 3, 2.5 \}$. The choice of $\Sigma_n$ and targets is detailed in the Appendix.

This experiment, which results are in Figure \ref{fig:nu}, strengthens the hypothesis of Ledoit and Wolf \cite{Ledoit2004} saying that bounded fourth moments is experimentally enough, and $\nu=3$ is what we found the be the frontier for the convergence to hold


\subsection{Application to non-stationary time series}

In order to give the practitioner a real example on an application, we consider a standard portfolio management problem in finance.

We use daily log-returns of $323$ assets that stayed in the S\&P500 between 2010 and 2022 \footnote{Data were retrieved from YahooFinance using \href{https://github.com/ranaroussi/yfinance}{yfinance} package, and available with the code on \href{https://github.com/cvcer/MTSE}{GitHub}.}. We recall that the log-return $r_t$ at day $t$ of a price $s_t$ is defined as the log-increment: $r_t = \log(s_t) - \log(s_{t-1})$. So we have a $323$-dimensional time series, with one vector per day.

Those $323$ assets are classified in $11$ industrial sectors, known as the \textbf{GICS sectors}. We expect to take advantage of it, using those sectors to build $11$ targets.

To evaluate the performance of different covariance estimators, we used the General Minimum Variance - GMV - portfolio problem, a standard method to compare covariance estimators. The idea of GMV is to determine at the beginning of each month $T$ a weight vector $w_T \in \R^{323}$, and invest $w_{i,T}$ on the asset $i$ for one month: this is our portfolio. 

At each day $t$ of the month $T$, around $20$ business days, we collect the log-returns of our portfolio $R_t = w_T^T r_t$, and at the end of the month we store the empirical variance of our portfolio during the month: $v_T = \V[\{R_t \}_{t=T}^{T+20}]$. The weight vector $w_T$ is designed to minimize this variance $v_T$, under constraint that $w_T^\top 1 = 1$. There is closed-form formula, where $\Sigma$ is the covariance matrix of the random vector $r_t$: $w_T =  \frac{\Sigma^{-1} 1}{1^\top \Sigma^{-1} 1}$. The pseudo-inverse is used when not invertible.

The experiment goes as following: at the beginning of each month $T$, we estimate $\hat \Sigma$ on the last $K$ months, \textit{i.e} roughly $20K$ data points, and we compute $\hat w_T= \frac{\hat \Sigma^{-1} 1}{1^T \hat \Sigma^{-1} 1}$. At the end of the month we store the empirical variance $\hat v_T$.  We do this procedure over $135$ months, between 2010 and 2022, and we compare the Cumulative Variance $\sum_{T=1}^{135} v_T$, the results are shown in the Table \ref{GMV}. The better the estimator $\hat \Sigma$, the lower is the Cumulative Variance.

The number $K$ of months used to fit the covariance estimators varies between $3$ and $15$. Along with sample covariance $S$, we compare:
\begin{itemize}
	\item MTSE, where the $11$ targets are chosen to be the Identity matrix on the sub-selection of assets belonging to the same GICS sector, and $0$ on the other indexes of the matrix,
	\item MTSE\_L defined in \cite{Lancewicki2014}, another multi-target estimator. The code being unavailable (as remarked in \cite{Gray2018}, Section 6), we implemented the algorithm specifically for our set of 11 GICS targets, using the fact that they are diagonal to get rid of the computation of $E_q \otimes E_q$ (Kronecker product) in \cite{Lancewicki2014} of size $327^2 \times 327^2$, which would be unusable in practice,
	\item MTSE\_G defined in \cite{Gray2018}, another multi-target estimator. In this setting, we cannot use our 11 GICS targets (as they are not invertible), so We used the nine targets proposed in their work,
	\item two STSE: LWO, Ledoit-Wolf estimator for unknown mean (see \cite{Oriol2023}), and OAS (see \cite{Chen2010}),
	\item two non-linear shrinkage: ANS (see \cite{Ledoit2020a}), and QIS (see \cite{Ledoit2019}).
\end{itemize}
The MTSE makes the difference in this difficult and non-stationnary task, as it consistently has a lower Cumulative Variance than the other considered algorithms.

\begin{table}
		
	\scalebox{0.90}{
	\begin{tabular}{|c|c|c|c|c|c|c|}
		\hline
		\text{Estimator}& $3$ \text{ months} & $4$ \text{ months} & $6$ \text{ months} &$ 9$ \text{ months} & $12$ \text{ months} & $15$ \text{ months}  \\
		\hline
		\text{MTSE} & $\mathbf{7.14 \times 10^{-3}}$ & $\mathbf{7.13 \times 10^{-3}}$ & $\mathbf{7.29 \times 10^{-3}}$& $\mathbf{7.09 \times 10^{-3}}$ & $\mathbf{7.62 \times 10^{-3}}$ & $\mathbf{7.52 \times 10^{-3}}$ \\
		\hline
		\text{MTSE\_L} & $7.71 \times 10^{-3}$ & $7.76 \times 10^{-3}$ & $8.41 \times 10^{-3}$& $8.28 \times 10^{-3}$ & $9.24 \times 10^{-3}$ & $8.66 \times 10^{-3}$ \\
		\hline
		\text{MTSE\_G} & $15.29 \times 10^{-3}$ & $15.28 \times 10^{-3}$ & $14.70 \times 10^{-3}$& $14.39 \times 10^{-3}$ & $14.30 \times 10^{-3}$ & $14.19 \times 10^{-3}$ \\
		\hline
		\text{LWO} & $7.36 \times 10^{-3}$ & $7.22 \times 10^{-3}$ & $7.40 \times 10^{-3}$& $7.25 \times 10^{-3}$ & $7.91 \times 10^{-3}$ & $7.65 \times 10^{-3}$ \\
		\hline
		\text{OAS}& $7.41 \times 10^{-3}$ & $7.34 \times 10^{-3}$ & $7.77 \times 10^{-3}$ & $ 7.90 \times 10^{-3}$ & $8.99 \times 10^{-3}$ & $8.90 \times 10^{-3}$ \\
		\hline
		\text{ANS} & $7.39 \times 10^{-3}$ & $7.34 \times 10^{-3}$ & $7.91 \times 10^{-3}$& $8.42 \times 10^{-3}$ & $10.95 \times 10^{-3}$ & $12.08 \times 10^{-3}$ \\
		\hline
		\text{QIS} & $1.38 \times 10^{0}$ & $33.22 \times 10^{-3}$ & $8.80 \times 10^{-3}$& $10.45 \times 10^{-3}$ & $17.56 \times 10^{-3}$ & $53.60 \times 10^{-3}$\\
		\hline
		\text{S} & $9.15 \times 10^{-3}$ & $8.64 \times 10^{-3}$ & $9.38 \times 10^{-3}$& $11.54 \times 10^{-3}$ & $16.40 \times 10^{-3}$ &$45.22 \times 10^{-3}$ \\
		\hline
		\end{tabular}
	}
	\caption{GMV Cumulative Variance for different estimators and number of months to fit them.}
	\label{GMV}
	\end{table}

\section{Limitations}
The strongest assumptions made in the theoretical analysis are that we have i.i.d. samples and bounded eighth moments. In practice, the i.i.d. assumption is hard to assess, particularly in time series analysis, as we saw in the financial dataset. The requirement of finite eighth moments can be a limitation, even if we saw that experimentally the convergence holds for heavier tails.

As in single-target estimation, the central point for performance in a task is a meaningful choice of targets. Some works address the topic in the literature, mostly for Single-Target estimators, see \cite{Ledoit2003, Wang2024} for example. For multi-target estimators, \cite{Lancewicki2014} proposes a set of four targets in Section III.B, and \cite{Gray2018} proposes a set of nine targets (section 2.4). However, a general method to design them and choose their number still lack in the literature to the best of our knowledge.

\section{Conclusion}
We derived the oracle estimator solving the multi-target linear shrinkage covariance estimator and proposed an explicit analytic estimator of it. The main result of the contribution is the proof of the convergence of the estimator under Kolmogorov asymptotics. 

The numerical simulations highlight the gain in performance to add targets to the estimation. Experimentally, we argue that the convergence holds with only bounded fourth moments. The experiment in a non-stationary setup underlines the adaptability of our estimator in those situations, compared to other estimators.

\section*{Statements and Declarations}
\begin{itemize}
	\item Funding: This research is funded by the CIFRE doctoral contract 2022/0536.
	\item Conflict of interest/Competing interests: No conflict of interests nor competing interests.
	\item Ethics approval and consent to participate: not applicable.
	\item Consent for publication: all co-authors consent for publication.
	\item Data availability: data is available online via yfinance API as explained in the corpus.
	\item Materials availability: not applicable.
	\item Code availability: code is available online via the GitHub linked in the corpus.
	\item Author contribution: the author contributed to every stage of the article.
	\end{itemize}

\begin{appendices}

	\section{Unknown mean formulas}
	We extend the \textit{bona fide} estimator to setting of unknown mean, exposed in the introduction.
	\begin{dft}[Unknown mean - Variance estimators]
		We suppose the mean to be \textbf{unknown}. We define:
		\begin{equation}
		\begin{aligned}
			&\bar{b}_n^2 = \frac{1}{n^2} \sum_{k=1}^n \left \lVert \frac{n}{n-1} \tilde{X}_{\cdot k}^n \left(\tilde{X}_{\cdot k}^n\right)^T - S_n \right \rVert_n^2, \\
			&\hat{\V}(S_n) = \frac{(n-1)^2}{(n-2)(n-3)} \bar{b}_n^2 - \frac{1}{n(n-2)} \lVert S_n \rVert_n^2 - \frac{n-1}{n(n-2)(n-3)} p_n\langle S_n, I_{p_n} \rangle_n^2.
		\end{aligned}
		\end{equation}
		Let $i \in \llbracket 1, N_n \rrbracket$. We define:
		\begin{equation}
		\begin{aligned}
			&\hat{\V}\left(\langle S_n, T_n^{(i)} \rangle_n\right) = q_n^{(0)} \bar{b}_{T^{(i)},n}^2 +q_n^{(1)}\left\langle S_nT_n^{(i)}, T_n^{(i)}S_n \right \rangle_n  - q_n^{(2)}\langle S_n, T_n^{(i)} \rangle_n^2,  \\
			&\text{with }\bar{b}_{T^{(i)},n}^2 = \frac{p_n}{(n-1)^2} \sum_{k=1}^n \left \langle \tilde{X}_{\cdot k}^n \left(\tilde{X}_{\cdot k}^n\right)^T, T_n^{(i)}\right\rangle_n^2 - \frac{1}{n}\left\langle S_nT_n^{(i)}, T_n^{(i)}S_n \right \rangle_n, \\
			& \text{and }q_n^{(0)} =  \frac{(n-1)^2}{p_n(n-2)(n-3)}, q_n^{(1)} =\frac{n-1}{p_nn(n-2)} , q_n^{(2)} =\frac{n^2 -2n -1}{n(n-2)(n-3)} .
		\end{aligned}
		\end{equation}
		\end{dft}
	We use these estimators directly in the formulas for $S_n^* =\mathcal{P}_{\mathcal{S}_n^+} \left( c_0 S_n + (1-c_0) \sum_{i=1}^{N_n} c_{i}T_n^{(i)} \right)$, within the coefficients:
	\begin{equation}
	\begin{aligned}
	&c_0 = \min\left((\hat{c}_0)_+, 1\right), c_{i} =  (1 - c_0)\langle S_n,  T_n^{(i)} \rangle_n,
	\end{aligned}
	\end{equation}
	\begin{equation*}
	\begin{aligned}
	\text{with }\hat{c}_{0} = &\frac{1}{d_n^2}\Bigg(\lVert S_n \rVert_n^2 - \hat{\V}(S_n) - \sum_{k=1}^{N_n} \left(\langle S_n, T_n^{(k)} \rangle_n^2 - \hat{\V} \left(\langle S_n, T_n^{(k)} \rangle_n\right) \right) \Bigg),\\
	d_n^2 =& \lVert S_n \rVert_n^2 - \sum_{k=1}^{N_n} \langle S_n, T_n^{(k)} \rangle_n^2.\\
	\end{aligned}
	\end{equation*}
	This simple change of definition for $\hat{\V}\left(\langle S_n, T_n^{(i)} \rangle_n\right)$ and $\hat{\V}(S_n)$, depending on the setting of known or unknown mean, conserves the convergence properties we stated in the theoretical part. Every result is proven for both settings.
	
	\section{Experimental settings and extensions}
	In order to improve the complexity of the algorithm, we used the following implementation trick.
	\begin{rmk}[Implementation trick]
		$\bar{b}_{T^{(i)},n}^2$ can be made fast to compute in Python with Numpy. We use a complex matrix square root of $T_n^{(i)}$, denoted $R_n^{(i)}$. We have a $O(p_n^3 + np_n)$-complexity formula:
		\begin{equation}
		\begin{aligned}
			&\sum_{k=1}^n \left \langle \tilde{X}_{\cdot k}^n \left(\tilde{X}_{\cdot k}^n\right)^T, T_n^{(i)}\right\rangle_n^2 = \frac{1}{p_n^2} \sum_{k=1}^n \left( \sum_{i=1}^{p_n} (R_n^{(i)}  \tilde{X})_{i,k}^2 \right)^2.
		\end{aligned}
		\end{equation}
		\end{rmk}
	Using this trick, the complexity of our algorithm aligns with the complexity of Single Target Linear Shrinkage, such as the Ledoit-Wolf estimator.
	
	\subsection{Experiment 1: Target set impact}
	$\Sigma$ is built as follows: this is a block diagonal matrix composed of 10 block matrices $(B_i)_{i=1}^{10}$ of shape $(5,5)$. For each $i \in \llbracket 1, 10 \rrbracket$, $B_i$ is drawn from a Wishart distribution with $5$ degrees of freedom and $\sqrt{11-i}I_5$ as scale matrix. 
	\begin{equation}\label{eq_block}
	\begin{aligned}
		\Sigma = \begin{pmatrix}
			B_1 &  & (0) \\
			 & \ddots &   \\
			(0) & & B_{10}
		\end{pmatrix}
	\end{aligned}
	\end{equation}
	The underlying distribution of $X$ is a multivariate t-distribution with $\nu = 9$.
	
	$\{T_1^{(i)}\}$ is built as following: the first target $T_1^{(1)} = I$, and for $i \in \llbracket 2,10 \rrbracket$, the target $T_1^{(i)}$ is block-diagonal; where $\Sigma$ has the block $B_i$, $T_1^{(i)}$ has an Identity block; elsewhere $T_1^{(i)}$ has null values. 
	\begin{equation} \label{eq_target}
	\begin{aligned}
		T_1^{(i)} = \begin{pmatrix}
			(0) &  \ldots & (0) \\
			\vdots& (\text{block }i:) I   & \vdots \\
			(0)  & \ldots &  (0)
		\end{pmatrix}
	\end{aligned}
	\end{equation}
	The set of target is orthogonalized with the Gram-Schmidt algorithm.
	
	For the second set of targets $(T_2)^{(i)}$, we apply a permutation $\sigma$ on both rows and columns of $T_1^{(i)}$, where $\sigma$ is defined as:
	\begin{equation}
	\begin{aligned}
		\sigma: j \in \llbracket 1, p \rrbracket \longrightarrow (j+2) \text{ mod } p.
	\end{aligned}
	\end{equation}
	For the target $i \in \llbracket 1, 10 \rrbracket$, at index $(k,l) \in \llbracket 1, p \rrbracket^2$, we define $(T_2^{(i)})_{\sigma(k),\sigma(l)} = (T_1^{(i)})_{k,l}$.
	
	
	\subsection{Experiment 3: Heavy-tail robustness}
	We consider $\Sigma$ block diagonal with $5$ blocks of same size drawn independently from Wishart distributions of $\lfloor p/5\rfloor$ degrees of freedom, and the set of targets is made of $5$ block diagonal matrices, target $i$ with the Identity as block $i$ and $0$ elsewhere. We compare, for $20000$ independent experiments using t-Students with  degrees of freedom. 
	
	\section{Theoretical proofs}
	In this section, we develop the proofs of the theoretical results of the work.
	
	\begin{rmk}
		As all the estimators in the following are invariant by translation in the case where we do not know the mean, we can suppose $\E[X] = 0$ to make the notation lighter. Indeed, we only work with $(\tilde X_n)_{ik} = (X_n)_{ik} - \frac{1}{n}\sum_{k'=1}^n (X_n)_{ik'}$, which is not affected by any change of mean.
	\end{rmk}
	\subsection{Proof of Proposition 1}
	The minimizer $c^* \in \R^{N_n+1}$ of:
	\begin{equation}
		\begin{aligned}
			\minimize_{\substack{c \in \R^{N_n+1}}} \left \lVert c_0 S_n + \sum_{k=1}^{N_n} c_{k} T_n^{(k)} - \Sigma_n \right \rVert_n^2.
		\end{aligned}
		\end{equation} 
	solves the following equation:
	\begin{equation}
	\begin{aligned}
		A_n c^* = b_n,
	\end{aligned}
	\end{equation}
	with:
	\begin{equation}
	\begin{aligned} &A_n = \setlength\arraycolsep{1pt}
		\begin{footnotesize}\begin{pmatrix} 
				\lVert S_n \rVert_n^2 & \langle S_n, T_n^{(1)} \rangle_n & \ldots &  \langle S_n, T_n^{(N_n)} \rangle \\ 
				\langle S_n, T_n^{(1)} \rangle_n & 1 &  & (0)\\
				\vdots & & \ddots & \\
				\langle S_n, T_n^{(N_n)} \rangle& (0) &  & 1
		\end{pmatrix}\end{footnotesize} ,\\
	&b_n = \left[\langle S_n, \Sigma_n \rangle_n, \langle T_n^{(1)}, \Sigma_n \rangle_n,...,\langle T_n^{(N_n)}, \Sigma_n \rangle_n \right]^T.
	\end{aligned}
	\end{equation}
	We have then:
	\begin{equation}
	\begin{aligned}
		&\det(A_n) = \lVert S_n \rVert_n^2 - \sum_{k=1}^{N_n} \langle S_n, T_n^{(k)} \rangle_n^2,\\
		&\text{com}(A_n)^T =  \setlength\arraycolsep{0pt}
		\begin{footnotesize}\begin{pmatrix}
				 1& - \langle S_n, T_n^{(1)} \rangle_n &\ldots &- \langle S_n, T_n^{(N)} \rangle_n & & \\ 
				 - \langle S_n, T_n^{(1)} \rangle_n  &\\
				 \vdots & &(B_n)&  \\
				 - \langle S_n, T_n^{(N)} \rangle_n  &  & &
		\end{pmatrix}\end{footnotesize},  \\
	\end{aligned}
	\end{equation}
	\begin{equation}
	\begin{aligned}
		\text{with } &B_n = \left(\lVert S_n \rVert_n^2 - \sum_{k=1}^{N_n} \langle S_n, T_n^{(k)} \rangle_n^2\right) I_{N_n} + QQ^T,\\
		&Q = \left[\langle S_n, T_n^{(i)} \rangle_n \right]_{i=1,...,N_n}^T.
	\end{aligned}
	\end{equation}
	Under Assumption 5, $A_n$ is invertible and we have: $A_n^{-1} = \frac{1}{\det(A_n)} \text{com}(A_n)^T$.
	
	\begin{rmk}
	In the particular case where $A_n$ is not invertible (not supposing Assumption 5), remove a target $T_n^{(i)}$ such as $\langle S_n, T_n^{(i)} \rangle \neq 0$ and the new built $A_n$ will be invertible.
	\end{rmk}
	
	We can conclude by a simple computation:
	\begin{equation}
	\begin{aligned}
	&c_0^* = \frac{1}{\det(A_n)}\left(\langle S_n, \Sigma_n \rangle - \sum_{k=1}^{N_n} \langle S_n, T_n^{(k)} \rangle\langle \Sigma_n, T_n^{(k)} \rangle \right),\\
	&c_i^* =  \langle \Sigma_n,  T_n^{(i)} \rangle - c_0^* \langle S_n,  T_n^{(i)} \rangle . \\
	\end{aligned}
	\end{equation}
	
	\subsection{Proof of Lemma 1}
	We prove firstly an essential lemma on the convergence of the coordinate estimators - of the form $\langle S_n, T_n^{(k)} \rangle_n$.
	\begin{lem}[Coordinate estimators] 
	Under Assumptions 1 to 3, for all $k \in \llbracket 1,N_n \rrbracket$, $\langle S_n, T_n^{(k)} \rangle_n$ is an unbiased estimator of $\langle \Sigma_n, T_n^{(k)} \rangle_n$. Moreover, it exists $M \in \R_+$, independent of $n$ such as $\sum_{k=1}^{N_n} \E\left[\langle S_n - \Sigma_n, T_n^{(k)} \rangle_n^4\right] \leq \frac{M}{n}K_1 K_2$, where $K_1$ and $K_2$ are defined in Assumptions 1 and 2.
	\end{lem}
	
	By linearity of the expectation, for all $k \in \llbracket 1,N_n \rrbracket$, $\langle S_n, T_n^{(k)} \rangle_n$ is un unbiased estimator of $\langle \Sigma_n, T_n^{(k)} \rangle_n$. 
	
	For the second part of the lemma, let $k \in  \llbracket 1,N_n \rrbracket$ and let denote $T = T_n^{(k)} $. The following proof relies upon the Proposition 1 and its method of graph construction from indices in the Preliminary combinatorial result in \cite{Oriol2023}. We recall the method and the proposition here.
	
	Let $K \in \N^*$, and $K$ indices $(k_1,...,k_K) \in \llbracket 1,n \rrbracket^K$. 
	
	Let's associate a graph with $K$ vertices $\mathcal{V} = \{ 1,...,K\}$ to this set of indices. The set of edges $\mathcal{E}$ is built as following: there is an edge between the node $a \in \mathcal{V}$ and $b \in \mathcal{V}$, $a \neq b$ (we don't allow self-loops), if the corresponding indices are equal, i.e if $k_a = k_b$. We finally define our graph $\mathcal{G} = (\mathcal{V}, \mathcal{E})$. 
	\begin{ppt}[Proposition 1 from \cite{Oriol2023}]
		Let $\mathcal{G}=(\mathcal{V},\mathcal{E})$ a graph with $K$ vertices generated from some indices $(k_1^{(0)}, ...,k_K^{(0)}) \in \llbracket 1, n \rrbracket^K$ with the procedure described previously. Suppose $\mathcal{G}$ has $C \in \llbracket 1,K \rrbracket$ connected components. Then, there are $\prod_{i=0}^{C-1}(n-i)$ set of indices $(k_1, ...,k_K) \in \llbracket 1, n \rrbracket^K$ which have the associated graph $\mathcal{G}$.
	\end{ppt}

	\paragraph{Known mean} When the mean is known, choosing $M' \in \R_+$ the number of graphs with $4$ nodes and $2$ connected components, we have:
	\begin{equation}
	\begin{aligned}
	\E\left[\langle S_n - \Sigma_n, T \rangle_n^4\right] =&\frac{1}{n^4} \sum_{k_1,k_2,k_3,k_4} \E\left[\prod_{m \in \{k_1,k_2,k_3,k_4\}}\left\langle X_{m}^{n}\left(X_{m}^{n}\right)^T - \Sigma_n, T \right\rangle_n \right] \\
	\E\left[\langle S_n - \Sigma_n, T \rangle_n^4\right] \leq &\frac{M'n(n-1)}{n^4} \E\left[\langle X_{1}^{n}\left(X_{1}^{n}\right)^T - \Sigma_n, T \rangle_n^2 \right]^2 \\
	&+ \frac{n}{n^4} \E\left[\langle X_{1}^{n}\left(X_{1}^{n}\right)^T - \Sigma_n, T \rangle_n^4 \right].
	\end{aligned}
	\end{equation}
	Moreover, we have:
	\begin{equation}
	\begin{aligned}
	\sum_{k=1}^{N_n} \E\left[\langle X_{1}^{n}\left(X_{1}^{n}\right)^T - \Sigma_n,  T_n^{(k)} \rangle_n^4 \right] & \leq \E\left[ \left(\sum_k\langle X_{1}^{n}\left(X_{1}^{n}\right)^T - \Sigma_n,  T_n^{(k)} \rangle_n^2\right)^2 \right]  \\
	& \leq \E\left[ \left\lVert X_{1}^{n}\left(X_{1}^{n}\right)^T - \Sigma_n \right\rVert_n^4 \right]  \\
	& \leq \E\left[ \left\lVert X_{1}^{n}\left(X_{1}^{n}\right)^T \right\rVert_n^4 \right]  \\ 
	& = \frac{1}{p_n^2}\E\left[ \left(\sum_i y_{i1}^2 \right)^4 \right]  \\
	& = p_n^2\E\left[ \left(\frac{1}{p_n}\sum_i y_{i1}^2 \right)^4 \right]  \\
	& \leq p_n^2 \frac{1}{p_n} \sum_i \E\left[y_{i1}^8 \right]  \\
	\sum_{k=1}^{N_n} \E\left[\langle X_{1}^{n}\left(X_{1}^{n}\right)^T - \Sigma_n,  T_n^{(k)} \rangle_n^4 \right] & \leq p_n^2 K_2.
	\end{aligned}
	\end{equation}
	And,
	\begin{equation}
	\begin{aligned}
	\sum_{k=1}^{N_n} \E\left[\langle X_{1}^{n}\left(X_{1}^{n}\right)^T,  T_n^{(k)} \rangle_n^2 \right]^2 & \leq \max_i \left(\E\left[\langle X_{1}^{n}\left(X_{1}^{n}\right)^T,  T_n^{(i)} \rangle_n^2 \right]\right)  \sum_{k=1}^{N_n} \E\left[\langle X_{1}^{n}\left(X_{1}^{n}\right)^T,  T_n^{(k)} \rangle_n^2 \right]  \\
	& \leq \max_i \left(\E\left[\langle X_{1}^{n}\left(X_{1}^{n}\right)^T,  T_n^{(i)} \rangle_n^2 \right]\right) \E\left[\lVert X_{1}^{n}\left(X_{1}^{n}\right)^T \rVert_n^2 \right]  \\
	\sum_{k=1}^{N_n} \E\left[\langle X_{1}^{n}\left(X_{1}^{n}\right)^T,  T_n^{(k)} \rangle_n^2 \right]^2 & \leq  \max_i \left(\E\left[\langle X_{1}^{n}\left(X_{1}^{n}\right)^T,  T_n^{(i)} \rangle_n^2 \right]\right) p_n \sqrt{K_2}.  \\
	\end{aligned}
	\end{equation}
	For any of the targets $T$, we have with Assumption 2: $\E\left[\langle X_{1}^{n}\left(X_{1}^{n}\right)^T, T \rangle_n^2\right] \leq \sqrt{K_2}$.
	Finally, using $M = M' + K_1$, we have:
	\begin{equation}
	\begin{aligned}
	\sum_{k=1}^{N_n} \E\left[\langle S_n - \Sigma_n, T_n^{(k)} \rangle_n^4\right] \leq \frac{M' + K_1}{n}K_1 K_2  = \frac{M}{n}K_1 K_2.
	\end{aligned}
	\end{equation}
	
	\paragraph{Unkown mean} When the mean is unknown, we have a larger decomposition. We use the same technics as in the known mean setting, with more terms and cross-products to handle: 
	\begin{equation}
	\begin{aligned}
	&\E\left[\langle S_n - \Sigma_n, T \rangle_n^4\right]  = \\
	&\frac{1}{n^4(n-1)^4} \sum_{\substack{k_1, k_2, k_3, k_4 \\ k_1', k_2', k_3', k_4'}} \E \left[ \prod_{s = 1}^4 \left\langle X_{k_s}^{n}\left(X_{k_s}^{n}\right)^T - \Sigma_n, T \right\rangle_n \right] && \text{Term (1)}\\
		& - 4 \E \left[\prod_{s = 1}^3 \left\langle X_{k_s}^{n}\left(X_{k_s}^{n}\right)^T - \Sigma_n, T \right\rangle_n \left\langle X_{k_4}^{n}\left(X_{k_4'}^{n}\right)^T - \delta_{k_4=k_4'}\Sigma_n, T \right\rangle_n  \right] && \text{Term (2)} \\
		& + 6 \E \left[\prod_{s = 1}^2 \left\langle X_{k_s}^{n}\left(X_{k_s}^{n}\right)^T - \Sigma_n, T \right\rangle_n \prod_{t = 3}^4\left\langle X_{k_t}^{n}\left(X_{k_t'}^{n}\right)^T - \delta_{k_t=k_t'}\Sigma_n, T \right\rangle_n \right] && \text{Term (3)}\\
		& - 4 \E \left[\left\langle X_{k_1}^{n}\left(X_{k_1}^{n}\right)^T - \Sigma_n, T \right\rangle_n  \prod_{t = 2}^4\left\langle X_{k_t}^{n}\left(X_{k_t'}^{n}\right)^T - \delta_{k_t=k_t'}\Sigma_n, T \right\rangle_n \right] && \text{Term (4)} \\
		& + \E \left[\prod_{t = 1}^4 \left\langle X_{k_t}^{n}\left(X_{k_t'}^{n}\right)^T - \delta_{k_t=k_t'}\Sigma_n, T \right\rangle_n\right].&& \text{Term (5)} \\
	\end{aligned}
	\end{equation}
	We prove the higher bound of the form $M''K_1K_2/n$ term by term, in absolute value, summing over all targets $(T_n^{(k)})_{k=1}^{N_n}$. For that, we decompose each term and bound them by expectations we can control.
	
	\paragraph{Term (1)} This term is equal to $\frac{n^4}{(n-1)^4}\E\left[\langle S_n^{(\text{known mean)}} - \Sigma_n, T \rangle_n^4\right] $. Supposing $n \geq 4$, we have $\frac{n^4}{(n-1)^4} \leq 4$, so, with $M_1'' = 4M$, we have:
	\begin{equation}
	\begin{aligned}
	& \sum_{k=1}^{N_n}\frac{1}{n^4(n-1)^4} \sum_{\substack{k_1, k_2, k_3, k_4 \\ k_1', k_2', k_3', k_4'}} \left|\E \left[ \prod_{s = 1}^4 \left\langle X_{k_s}^{n}\left(X_{k_s}^{n}\right)^T - \Sigma_n, T \right\rangle_n \right]\right|  \leq \frac{M''_1}{n} K_1 K_2.
	\end{aligned}
	\end{equation}
	
	\paragraph{Term (2)} For this term, we have the following development, when discarding $0$-expectation terms:
	\begin{equation}
		\footnotesize
	\begin{aligned}
	& &&\frac{4}{n^4(n-1)^4} \sum_{\substack{k_1, k_2, k_3, k_4 \\ k_1', k_2', k_3', k_4'}}\left|\E \left[\prod_{s = 1}^3 \left\langle X_{k_s}^{n}\left(X_{k_s}^{n}\right)^T - \Sigma_n, T \right\rangle_n \left\langle X_{k_4}^{n}\left(X_{k_4'}^{n}\right)^T - \delta_{k_4=k_4'}\Sigma_n, T \right\rangle_n  \right]\right| \\
	&\leq && \frac{4}{n} \times \left|\text{Term(1)}\right| + \frac{24}{(n-1)^3}\left|\E\left[\left\langle X_{1}^{n}\left(X_{1}^{n}\right)^T - \Sigma_n, T \right \rangle_n^2 \left\langle X_{2}^{n}\left(X_{2}^{n}\right)^T - \Sigma_n, T \right \rangle_n \left\langle X_{1}^{n}\left(X_{2}^{n}\right)^T, T \right \rangle_n \right] \right| \\
	&\leq &&\frac{4}{n} \times \left|\text{Term(1)}\right| + \frac{24}{(n-1)^3}\sqrt{\E\left[\left\langle X_{1}^{n}\left(X_{1}^{n}\right)^T - \Sigma_n, T \right \rangle_n^4 \right]^{3/2}\E\left[\left\langle X_{1}^{n}\left(X_{2}^{n}\right)^T, T \right \rangle_n^2 \right]} \\
	&\leq &&\frac{4}{n} \times \left|\text{Term(1)}\right| + \frac{6}{(n-1)^3}\left(3\E\left[\left\langle X_{1}^{n}\left(X_{1}^{n}\right)^T - \Sigma_n, T \right \rangle_n^4 \right] + \E\left[\left\langle X_{1}^{n}\left(X_{2}^{n}\right)^T, T \right \rangle_n^2 \right]^2 \right).
	\end{aligned}
	\end{equation}
	At this point, we recall that we have to sum over each target $(T_n^{(k)})_{k=1}^{N_n}$ and then bound by a term of the form $M''_2K_1K_2/n$ for some $M_2'' \in \R$. $\text{Term(1)}$ is obviously controlled by the previous part. For the rest, we use the two following inequalities:
	\begin{itemize}
		\item from the proof in the known mean setting, we have that:
			\begin{equation}
			\begin{aligned}
			\sum_{k=1}^{N_n} \E\left[\left\langle X_{1}^{n}\left(X_{1}^{n}\right)^T - \Sigma_n, T_n^{(k)} \right \rangle_n^4 \right] \leq p_n^2 K_2,
			\end{aligned}
			\end{equation}
		\item and:
			\begin{equation}
			\begin{aligned}
			&\sum_{k=1}^{N_n} \E\left[\langle X_{1}^{n}\left(X_{2}^{n}\right)^T,  T_n^{(k)} \rangle_n^2 \right]^2 \\
			& \leq  \max_i \left(\E\left[\langle X_{1}^{n}\left(X_{2}^{n}\right)^T,  T_n^{(i)} \rangle_n^2 \right]\right) \sum_{k=1}^{N_n} \E\left[\langle X_{1}^{n}\left(X_{2}^{n}\right)^T,  T_n^{(k)} \rangle_n^2 \right]  \\
			& \leq  \max_i \left(\E\left[\frac{1}{p_n^2}\left(X_{2}^{n}\right)^T  T_n^{(i)} X_{1}^{n} \left(X_{1}^{n}\right)^T T_n^{(i)} X_{2}^{n} \right]\right) \E\left[\lVert X_{1}^{n}\left(X_{2}^{n}\right)^T \rVert_n^2 \right]  \\
			 & \leq  \max_i \left(\frac{1}{p_n} \left\langle \Sigma_n  T_n^{(i)},  T_n^{(i)} \Sigma_n\right\rangle_n \right) p_n \sqrt{K_2}  \\
			  & \leq  \max_i \left(\frac{1}{p_n} \left\lVert \Sigma_n  T_n^{(i)}\right\rVert_n^2 \right) p_n \sqrt{K_2}  \\
			  & \leq \left\lVert \Sigma_n \right\rVert_n^2 p_n \sqrt{K_2}  \\
			&\leq  p_n K_2.
			\end{aligned}
			\end{equation}
	\end{itemize}
	Supposing $n \geq 4$, with $M_2'' = M_1'' + 54 K_1$, we then have:
	\begin{equation}
	\begin{aligned}
	& \sum_{k=1}^{N_n}\frac{4}{n^4(n-1)^4} \sum_{\substack{k_1, k_2, k_3, k_4 \\ k_1', k_2', k_3', k_4'}}\left|\E \left[\prod_{s = 1}^3 \left\langle X_{k_s}^{n}\left(X_{k_s}^{n}\right)^T - \Sigma_n, T \right\rangle_n \left\langle X_{k_4}^{n}\left(X_{k_4'}^{n}\right)^T - \delta_{k_4=k_4'}\Sigma_n, T \right\rangle_n  \right]\right|  \\
	&\leq \frac{M''_2}{n} K_1 K_2.
	\end{aligned}
	\end{equation}

	\paragraph{Term (3)} Similarly for this term, we obtain:
	\begin{equation}
		\footnotesize
	\begin{aligned}
	& &&\frac{6}{n^4(n-1)^4} \sum_{\substack{k_1, k_2, k_3, k_4 \\ k_1', k_2', k_3', k_4'}}\left|\E \left[\prod_{s = 1}^2 \left\langle X_{k_s}^{n}\left(X_{k_s}^{n}\right)^T - \Sigma_n, T \right\rangle_n \prod_{t = 3}^4\left\langle X_{k_t}^{n}\left(X_{k_t'}^{n}\right)^T - \delta_{k_t=k_t'}\Sigma_n, T \right\rangle_n \right]\right| \\
	&\leq &&\frac{6(2n-1)}{n^2} \times \left|\text{Term(2)}\right|  \\
	& &&+ \frac{12(n-2)}{n(n-1)^3}\left|\E\left[\left\langle X_{1}^{n}\left(X_{1}^{n}\right)^T - \Sigma_n, T \right \rangle_n^2 \left\langle X_{2}^{n}\left(X_{3}^{n}\right)^T - \Sigma_n, T \right \rangle_n^2  \right]\right|   \\
	& &&+ \frac{48(n-2)}{n(n-1)^3}\left|\E\left[\left\langle X_{1}^{n}\left(X_{1}^{n}\right)^T - \Sigma_n, T \right \rangle_n \left\langle X_{2}^{n}\left(X_{2}^{n}\right)^T - \Sigma_n, T \right \rangle_n \left \langle X_{1}^{n}\left(X_{3}^{n}\right)^T, T \right \rangle_n \left\langle X_{2}^{n}\left(X_{3}^{n}\right)^T, T \right \rangle_n  \right]\right|  \\
	& &&+ \frac{24}{n(n-1)^3}\left|\E\left[\left\langle X_{1}^{n}\left(X_{1}^{n}\right)^T - \Sigma_n, T \right \rangle_n \left\langle X_{2}^{n}\left(X_{2}^{n}\right)^T - \Sigma_n, T \right \rangle_n \left \langle X_{1}^{n}\left(X_{2}^{n}\right)^T, T \right \rangle_n^2  \right]\right|  \\
	& &&+ \frac{24}{n(n-1)^3}\left||\E\left[\left\langle X_{1}^{n}\left(X_{1}^{n}\right)^T - \Sigma_n, T \right \rangle_n^2 \left \langle X_{1}^{n}\left(X_{2}^{n}\right)^T, T \right \rangle_n^2  \right]\right| \\
	& \leq && \frac{6(2n-1)}{n^2} \times \left|\text{Term(2)}\right| 
	+ \frac{24}{n(n-1)^2}\left(\E\left[\left\langle X_{1}^{n}\left(X_{1}^{n}\right)^T - \Sigma_n, T \right \rangle_n^4 \right] 
	+ \E\left[\left\langle X_{1}^{n}\left(X_{2}^{n}\right)^T, T \right \rangle_n^4 \right] \right)\\
	& &&+ \frac{6(n-2)}{n(n-1)^3}\left( \E\left[\left\langle X_{1}^{n}\left(X_{1}^{n}\right)^T - \Sigma_n, T \right \rangle_n^2 \right]^2
	+ \E\left[\left\langle X_{1}^{n}\left(X_{2}^{n}\right)^T, T \right \rangle_n^2 \right]^2 \right).
	\end{aligned}
	\end{equation}
	We use the two following inequalities to finish the work, combined with the two ones we used for Term (2):
	\begin{itemize}
		\item from the proof in the known mean setting, we have that:
			\begin{equation}
			\begin{aligned}
			\sum_{k=1}^{N_n} \E\left[\left\langle X_{1}^{n}\left(X_{1}^{n}\right)^T - \Sigma_n, T_n^{(k)} \right \rangle_n^2 \right]^2 \leq p_n K_2,
			\end{aligned}
			\end{equation}
		\item and:
			\begin{equation}
			\begin{aligned}
			\sum_{k=1}^{N_n} \E\left[\langle X_{1}^{n}\left(X_{2}^{n}\right)^T,  T_n^{(k)} \rangle_n^4 \right] & \leq \E\left[ \left(\sum_{k=1}^{N_n} \langle X_{1}^{n}\left(X_{2}^{n}\right)^T,  T_n^{(k)} \rangle_n^2\right)^2 \right]  \\
			& \leq \E\left[ \left\lVert X_{1}^{n}\left(X_{2}^{n}\right)^T \right\rVert_n^4 \right]  \\
			& = \frac{1}{p_n^2}\E\left[ \left(\sum_i y_{i1}y_{i2} \right)^4 \right]  \\
			& \leq p_n^2 \frac{1}{p_n} \sum_i \E\left[y_{i1}^4y_{i2}^4 \right]  \\
			& \leq p_n^2 \frac{1}{p_n} \sum_i \E\left[y_{i1}^8 \right]  \\
			\sum_{k=1}^{N_n} \E\left[\langle X_{1}^{n}\left(X_{2}^{n}\right)^T,  T_n^{(k)} \rangle_n^4 \right]& \leq p_n^2 K_2.
			\end{aligned}
			\end{equation}
	\end{itemize}
	Supposing $n \geq 4$, with $M_3'' = 3M_2'' + 86 K_1 + 1$, we then have:
	\begin{equation}
		\small
	\begin{aligned}
	& \sum_{k=1}^{N_n} \frac{6}{n^4(n-1)^4} \sum_{\substack{k_1, k_2, k_3, k_4 \\ k_1', k_2', k_3', k_4'}}\left|\E \left[\prod_{s = 1}^2 \left\langle X_{k_s}^{n}\left(X_{k_s}^{n}\right)^T - \Sigma_n, T \right\rangle_n \prod_{t = 3}^4\left\langle X_{k_t}^{n}\left(X_{k_t'}^{n}\right)^T - \delta_{k_t=k_t'}\Sigma_n, T \right\rangle_n \right]\right| \\
	&\leq \frac{M''_3}{n} K_1 K_2.
	\end{aligned}
	\end{equation}

	\paragraph{Term (4)} For that one, we have:
	\begin{equation}
		\footnotesize
	\begin{aligned}
	& &&\frac{4}{n^4(n-1)^4} \sum_{\substack{k_1, k_2, k_3, k_4 \\ k_1', k_2', k_3', k_4'}}\left|\E \left[\left\langle X_{k_1}^{n}\left(X_{k_1}^{n}\right)^T - \Sigma_n, T \right\rangle_n  \prod_{t = 2}^4\left\langle X_{k_t}^{n}\left(X_{k_t'}^{n}\right)^T - \delta_{k_t=k_t'}\Sigma_n, T \right\rangle_n \right]\right| \\
	& \leq &&\frac{4(n-1)(n+3)}{n^3} \times \left|\text{Term(3)}\right| \\
	& && + \frac{96(n-2)}{n^2(n-1)^3} \left|\E\left[\left\langle X_{1}^{n}\left(X_{1}^{n}\right)^T - \Sigma_n, T \right \rangle_n \left\langle X_{1}^{n}\left(X_{2}^{n}\right)^T, T \right \rangle_n \left \langle X_{2}^{n}\left(X_{3}^{n}\right)^T, T \right \rangle_n \left\langle X_{3}^{n}\left(X_{1}^{n}\right)^T, T \right \rangle_n  \right]\right| \\
	& && + \frac{96(n-2)}{n^2(n-1)^3} \left|\E\left[\left\langle X_{1}^{n}\left(X_{1}^{n}\right)^T - \Sigma_n, T \right \rangle_n \left\langle X_{1}^{n}\left(X_{2}^{n}\right)^T, T \right \rangle_n \left \langle X_{2}^{n}\left(X_{3}^{n}\right)^T, T \right \rangle_n^2\right]\right| \\
	& && + \frac{24}{n^2(n-1)^3} \left|\E\left[\left\langle X_{1}^{n}\left(X_{1}^{n}\right)^T - \Sigma_n, T \right \rangle_n \left\langle X_{1}^{n}\left(X_{2}^{n}\right)^T, T \right \rangle_n^3\right]\right| \\
	& \leq && \frac{4(n-1)(n+3)}{n^3} \times \left|\text{Term(3)}\right| + \\
	& &&\frac{48n - 88}{n^2(n-1)^3} \left(3\E\left[\left\langle X_{1}^{n}\left(X_{2}^{n}\right)^T, T \right \rangle_n^4 \right]+ \E\left[\left\langle X_{1}^{n}\left(X_{1}^{n}\right)^T - \Sigma_n, T \right \rangle_n^2 \right]^2 \right)
	\end{aligned}
	\end{equation}
	Supposing $n \geq 4$, with $M_4'' = M_3'' + 47 K_1 + 4 $, we then have:
	\begin{equation}
		\footnotesize
	\begin{aligned}
	& \sum_{k=1}^{N_n} \frac{4}{n^4(n-1)^4} \sum_{\substack{k_1, k_2, k_3, k_4 \\ k_1', k_2', k_3', k_4'}}\left|\E \left[\left\langle X_{k_1}^{n}\left(X_{k_1}^{n}\right)^T - \Sigma_n, T_k^{(n)} \right\rangle_n  \prod_{t = 2}^4\left\langle X_{k_t}^{n}\left(X_{k_t'}^{n}\right)^T - \delta_{k_t=k_t'}\Sigma_n, T_k^{(n)} \right\rangle_n \right]\right| \\
	&\leq \frac{M''_4}{n} K_1 K_2.
	\end{aligned}
	\end{equation}
	
	\paragraph{Term (5)} Following the same idea, we have for the last term:
	\begin{equation}
		\footnotesize
	\begin{aligned}
	& &&\frac{1}{n^4(n-1)^4} \sum_{\substack{k_1, k_2, k_3, k_4 \\ k_1', k_2', k_3', k_4'}}\left|\E \left[\prod_{t = 1}^4 \left\langle X_{k_t}^{n}\left(X_{k_t'}^{n}\right)^T - \delta_{k_t=k_t'}\Sigma_n, T \right\rangle_n\right]\right| \\
	& \leq && \frac{n^3+n^2-4}{n^4} \times \left|\text{Term(4)}\right| 
	+ \frac{8}{n^3(n-1)^3} \left|\E\left[\left\langle X_{1}^{n}\left(X_{2}^{n}\right)^T, T \right \rangle_n^4\right]\right|  \\
	& &&+ \frac{24(n-2)}{n^3(n-1)^3} \left|\E\left[\left\langle X_{1}^{n}\left(X_{2}^{n}\right)^T, T \right \rangle_n^2 \left\langle X_{1}^{n}\left(X_{3}^{n}\right)^T, T \right \rangle_n^2\right]\right| \\
	& &&+ \frac{12(n-2)(n-3)}{n^3(n-1)^3} \left|\E\left[\left\langle X_{1}^{n}\left(X_{2}^{n}\right)^T, T \right \rangle_n^2 \left\langle X_{3}^{n}\left(X_{4}^{n}\right)^T, T \right \rangle_n^2\right]\right| \\
	& &&+ \frac{48(n-2)(n-3)}{n^3(n-1)^3} \left|\E\left[\left\langle X_{1}^{n}\left(X_{2}^{n}\right)^T, T \right \rangle_n^2 \left\langle X_{2}^{n}\left(X_{3}^{n}\right)^T, T \right \rangle_n \left\langle X_{3}^{n}\left(X_{1}^{n}\right)^T, T \right \rangle_n\right]\right| \\
	& &&+ \frac{48(n-2)(n-3)}{n^3(n-1)^3} \left|\E\left[\left\langle X_{1}^{n}\left(X_{2}^{n}\right)^T, T \right \rangle_n \left\langle X_{2}^{n}\left(X_{3}^{n}\right)^T, T \right \rangle_n \left\langle X_{3}^{n}\left(X_{4}^{n}\right)^T, T \right \rangle_n \left\langle X_{4}^{n}\left(X_{1}^{n}\right)^T, T \right \rangle_n\right]\right| \\
	& \leq &&  \frac{n^3+n^2-4}{n^4} \times \left|\text{Term(4)}\right| + \frac{72n+8}{n^3(n-1)^3} \E\left[\left\langle X_{1}^{n}\left(X_{2}^{n}\right)^T, T \right \rangle_n^4 \right] \\
	& &&+ \frac{60(n-2)(n-3)}{n^3(n-1)^3} \E\left[\left\langle X_{1}^{n}\left(X_{2}^{n}\right)^T, T \right \rangle_n^2 \right]^2.
	\end{aligned}
	\end{equation}
	Supposing $n \geq 4$, with $M_5'' = M_4'' + 11 K_1 + 2$, we then have:
	\begin{equation}
	\begin{aligned}\label{eq_m5}
	& \sum_{k=1}^{N_n} \frac{1}{n^4(n-1)^4} \sum_{\substack{k_1, k_2, k_3, k_4 \\ k_1', k_2', k_3', k_4'}}\left|\E \left[\prod_{t = 1}^4 \left\langle X_{k_t}^{n}\left(X_{k_t'}^{n}\right)^T - \delta_{k_t=k_t'}\Sigma_n, T_k^{(n)} \right\rangle_n\right]\right| \leq \frac{M''_5}{n} K_1 K_2.
	\end{aligned}
	\end{equation}
	
	Finally, combing the 5 terms, we have with $M'' = \sum_{i=1}^5 M_i''$:
	\begin{equation}
	\begin{aligned}
	\sum_{k=1}^{N_n} \E\left[\langle S_n - \Sigma_n, T_n^{(k)} \rangle_n^4\right] \leq \frac{M''}{n}K_1 K_2,
	\end{aligned}
	\end{equation}
	which concludes the proof for the unknown mean.
	
	\subsection{Proof of Theorem 1}
	\paragraph{Known mean} When the mean is known, $\hat{\V}(S_n)$ is trivially unbiased when developping. $\E\left[\left(\hat{\V}(S_n) - \V(S_n)\right)^2\right] \longrightarrow 0$ is a consequence of Lemma 3.4 in \cite{Ledoit2004} about the quadratic mean convergence of $\bar b_n^2 - \beta_n^2 = \frac{n-1}{n}\hat{\V}(S_n) - \V(S_n)$.
	
	$\hat{\V}(\langle S_n, T_n^{(i)} \rangle_n)$ is trivially unbiased when developing. Let's prove that $ \E\left[\left(\sum_i \hat{\V}(\langle S_n, T_n^{(i)} \rangle_n) - \V(\langle S_n, T_n^{(i)} \rangle_n)\right)^2\right] \longrightarrow 0$.
	\begin{equation}
	\begin{aligned}
	&\E\left[\left(\sum_i \frac{1}{n(n-1)} \sum_k \left \langle X_{k}^{n} (X_{k}^{n})^T - \Sigma_n,T_n^{(i)} \right \rangle_n^2 - \E\left[\left \langle X_{k}^{n} (X_{k}^{n})^T - \Sigma_n,T_n^{(i)} \right \rangle_n^2 \right] \right)^2 \right] \\
	&= \frac{1}{n(n-1)^2} \E\left[\left(\sum_i \left \langle X_{1}^{n} (X_{1}^{n})^T - \Sigma_n,T_n^{(i)} \right \rangle_n^2 - \E\left[\left \langle X_{1}^{n} (X_{1}^{n})^T - \Sigma_n,T_n^{(i)} \right \rangle_n^2 \right] \right)^2 \right]\\
	&\leq \frac{N_n}{n(n-1)^2} \E\left[\sum_i \left \langle X_{1}^{n} (X_{1}^{n})^T - \Sigma_n,T_n^{(i)} \right \rangle_n^4 \right]  \\
	&\leq \frac{N_n}{n(n-1)^2} \times 7 \E\left[\sum_i \left \langle X_{1}^{n} (X_{1}^{n})^T,T_n^{(i)} \right \rangle_n^4 \right]  \\
	&\leq 7\frac{nN_n}{(n-1)^2} K_1^2 K_2.
	\end{aligned}
	\end{equation}
	And,
	\begin{equation}
	\begin{aligned}
	&\E\left[\left(\sum_i \frac{1}{n-1}\left \langle S_n - \Sigma_n,T_n^{(i)} \right \rangle_n^2 - \E\left[\frac{1}{n-1}\left \langle S_n - \Sigma_n,T_n^{(i)} \right \rangle_n^2 \right] \right)^2 \right] \\
	&\leq \frac{N_n}{(n-1)^2}\sum_i \E\left[\left \langle S_n - \Sigma_n,T_n^{(i)} \right \rangle_n^4 \right]\\
	&\leq \frac{N_nMK_1K_2}{n(n-1)^2}. &&\text{(Lemma 1)}
	\end{aligned}
	\end{equation}
	Which concludes the proof when the mean is known.
	
	\paragraph{Unknown mean} 
	When the mean is unknown, $\hat{\V}(S_n)$ is unbiased thanks to Lemma 9 in \cite{Oriol2023}. $\E\left[\left(\hat{\V}(S_n) - \V(S_n)\right)^2\right] \longrightarrow 0$ is a consequence of Lemma 9 in \cite{Oriol2023} about the quadratic mean convergence of $b_n^2 - \beta_n^2 = \hat{\V}(S_n) - \V(S_n)$.
	
	$\hat{\V}(\langle S_n, T_n^{(i)} \rangle_n)$ is unbiased thanks to Lemma 8 in \cite{Oriol2023} and the development of the variance. Let's prove that $ \E\left[\left(\sum_i \hat{\V}(\langle S_n, T_n^{(i)} \rangle_n) - \V(\langle S_n, T_n^{(i)} \rangle_n)\right)^2\right] \longrightarrow 0$. 
	For that, we need to prove the convergence of three terms:
	\begin{itemize}
		\item[(I.)] $\V\left[\frac{1}{(n-1)^2}\sum_{i=1}^{N_n} \sum_{k=1}^n \langle \tilde{X}_{k}^{n} (\tilde{X}_{k}^{n})^T, T_n^{(i)} \rangle_n^2 \right] \longrightarrow 0$,
		\item[(II.)] $\V\left[\frac{1}{pn}\sum_{i=1}^{N_n} \langle S_n T_n^{(i)}, T_n^{(i)} S_n \rangle_n \right] \longrightarrow 0$,
		\item[(III.)] $\V\left[\frac{1}{n}\sum_{i=1}^{N_n} \langle S_n, T_n^{(i)} \rangle_n^2 \right] \longrightarrow 0$.
	\end{itemize}
	
	\paragraph{Term (I)} We first show that $\V\left[\frac{1}{(n-1)^2}\sum_{i=1}^{N_n} \sum_{k=1}^n \langle \tilde{X}_{k}^{n} (\tilde{X}_{k}^{n})^T, T_n^{(i)} \rangle_n^2 \right] \longrightarrow 0$.
	\begin{equation}
	\begin{aligned}
	& &&\V\left[\frac{1}{(n-1)^2}\sum_{i=1}^{N_n} \sum_{k=1}^n \langle \tilde{X}_{k}^{n} (\tilde{X}_{k}^{n})^T, T_n^{(i)} \rangle_n^2 \right] \\
	&\leq && \E\left[\left(\frac{1}{(n-1)^2}\sum_{i=1}^{N_n} \sum_{k=1}^n \langle \tilde{X}_{k}^{n} (\tilde{X}_{k}^{n})^T, T_n^{(i)} \rangle_n^2 \right)^2\right] \\
	&\leq && \frac{N_n}{(n-1)^4} \sum_{i=1}^{N_n}\E\left[\left(\sum_{k=1}^n \langle \tilde{X}_{k}^{n} (\tilde{X}_{k}^{n})^T, T_n^{(i)} \rangle_n^2 \right)^2\right] \\
	&= &&\frac{N_n}{(n-1)^4} \sum_{i=1}^{N_n} \E\Bigg[\Bigg(\sum_{k=1}^n 
	\langle X_{k}^{n} (X_{k}^{n})^T, T_n^{(i)} \rangle_n^2 
	-4 \langle X_{k}^{n} (X_{k}^{n})^T, T_n^{(i)} \rangle_n\langle X_{k}^{n} (\bar{X}^{n})^T, T_n^{(i)} \rangle_n \\
	& &&+4\langle X_{k}^{n} (\bar{X}^{n})^T, T_n^{(i)} \rangle_n^2
	+\langle \bar{X}^{n} (\bar{X}^{n})^T, T_n^{(i)} \rangle_n^2
	-4\langle X_{k}^{n} (\bar{X}^{n})^T, T_n^{(i)} \rangle_n \langle \bar{X}^{n} (\bar{X}^{n})^T, T_n^{(i)} \rangle_n \\
	& &&+2 \langle X_{k}^{n} (X_{k}^{n})^T, T_n^{(i)} \rangle_n \langle \bar{X}^{n} (\bar{X}^{n})^T, T_n^{(i)} \rangle_n
	\Bigg)^2\Bigg] \\
	\end{aligned}
	\end{equation}
	So, we have:
	\begin{equation}
	\begin{aligned}
	& &&\V\left[\frac{1}{(n-1)^2}\sum_{i=1}^{N_n} \sum_{k=1}^n \langle \tilde{X}_{k}^{n} (\tilde{X}_{k}^{n})^T, T_n^{(i)} \rangle_n^2 \right] \\
	& \leq &&\frac{6N_n}{(n-1)^4} \sum_{i=1}^{N_n}\E\left[\left(\sum_{k=1}^n \langle X_{k}^{n} (X_{k}^{n})^T, T_n^{(i)} \rangle_n^2 \right)^2\right] &&\text{Term (I.1)}\\
	& &&+ \frac{6N_n}{(n-1)^4} \sum_{i=1}^{N_n}\E\left[\left(\sum_{k=1}^n 4 \langle X_{k}^{n} (X_{k}^{n})^T, T_n^{(i)} \rangle_n\langle X_{k}^{n} (\bar{X}^{n})^T, T_n^{(i)} \rangle_n\right)^2\right] &&\text{Term (I.2)}\\
	& &&+\frac{6N_n}{(n-1)^4} \sum_{i=1}^{N_n}\E\left[\left(\sum_{k=1}^n 4\langle X_{k}^{n} (\bar{X}^{n})^T, T_n^{(i)} \rangle_n^2\right)^2\right] &&\text{Term (I.3)}\\
	& &&+\frac{6N_n}{(n-1)^4} \sum_{i=1}^{N_n}\E\left[\left(\sum_{k=1}^n\langle \bar{X}^{n} (\bar{X}^{n})^T, T_n^{(i)} \rangle_n^2\right)^2\right] &&\text{Term (I.4)}\\
	& &&+\frac{6N_n}{(n-1)^4} \sum_{i=1}^{N_n}\E\left[ \left(\sum_{k=1}^n 4\langle X_{k}^{n} (\bar{X}^{n})^T, T_n^{(i)} \rangle_n \langle \bar{X}^{n} (\bar{X}^{n})^T, T_n^{(i)} \rangle_n\right)^2\right] &&\text{Term (I.5)}\\
	& &&+ \frac{6N_n}{(n-1)^4} \sum_{i=1}^{N_n}\E\left[\left(\sum_{k=1}^n2 \langle X_{k}^{n} (X_{k}^{n})^T, T_n^{(i)} \rangle_n \langle \bar{X}^{n} (\bar{X}^{n})^T, T_n^{(i)} \rangle_n\right)^2\right]&&\text{Term (I.6)}.
	\end{aligned}
	\end{equation}
	With this decomposition, we have $6$ distinct terms, and we prove the convergence to $0$ of each one.
	
	\paragraph{Term (I.1)}
	\begin{equation}
	\begin{aligned}
	&\frac{6N_n}{(n-1)^4} \sum_{i=1}^{N_n}\E\left[\left(\sum_{k=1}^n \langle X_{k}^{n} (X_{k}^{n})^T, T_n^{(i)} \rangle_n^2 \right)^2\right] 
	&& \leq \frac{6nN_n}{(n-1)^4}  \sum_{i=1}^{N_n}\E\left[\langle X_{1}^{n} (X_{1}^{n})^T, T_n^{(i)} \rangle_n^4\right] \\
	& &&\leq \frac{6nN_n p_n^2 K_2}{(n-1)^2} \\
	&\frac{6N_n}{(n-1)^4} \sum_{i=1}^{N_n}\E\left[\left(\sum_{k=1}^n \langle X_{k}^{n} (X_{k}^{n})^T, T_n^{(i)} \rangle_n^2 \right)^2\right] && \longrightarrow 0.
	\end{aligned}
	\end{equation}
	
	\paragraph{Term (I.2)}
	\begin{equation}
	\begin{aligned}
	& &&\frac{6N_n}{(n-1)^4} \sum_{i=1}^{N_n}\E\left[\left(\sum_{k=1}^n 4 \langle X_{k}^{n} (X_{k}^{n})^T, T_n^{(i)} \rangle_n\langle X_{k}^{n} (\bar{X}^{n})^T, T_n^{(i)} \rangle_n\right)^2\right] \\
	&= &&\frac{6\times 16N_n}{n(n-1)^4} \sum_{i=1}^{N_n}\E\left[\langle X_{1}^{n} (X_{1}^{n})^T, T_n^{(i)} \rangle_n^4\right] \\
	& &&+ \frac{6\times 16N_n}{n(n-1)^3} \sum_{i=1}^{N_n}\E\left[\langle X_{1}^{n} (X_{1}^{n})^T, T_n^{(i)} \rangle_n^2\langle X_{1}^{n} (X_{2}^{n})^T, T_n^{(i)} \rangle_n^2\right] \\
	&\leq &&\frac{96 n N_n}{(n-1)^4} K_1^2 K_2 + \frac{96 n N_n}{(n-1)^3} K_1^2 K_2  \longrightarrow 0.
	\end{aligned}
	\end{equation}
	
	\paragraph{Term (I.3)}
	\begin{equation}
	\begin{aligned}
	& &&\frac{6N_n}{(n-1)^4} \sum_{i=1}^{N_n}\E\left[\left(\sum_{k=1}^n 4\langle X_{k}^{n} (\bar{X}^{n})^T, T_n^{(i)} \rangle_n^2\right)^2\right] \\
	&\leq && \frac{96n^2N_n}{(n-1)^4} \sum_{i=1}^{N_n}\E\left[\langle X_{1}^{n} (\bar{X}^{n})^T, T_n^{(i)} \rangle_n^4\right] \\
	&= &&\frac{96N_n}{n^2(n-1)^4} \sum_{i=1}^{N_n} \E\Bigg[
	\langle X_{1}^{n} (X_{1}^{n})^T, T_n^{(i)} \rangle_n^4\\
	& &&+6(n-1)\langle X_{1}^{n} (X_{1}^{n})^T, T_n^{(i)} \rangle_n^2\langle X_{1}^{n} (X_{2}^{n})^T, T_n^{(i)} \rangle_n^2 \\
	& &&
	+3(n-1)(n-2)\langle X_{1}^{n} (X_{2}^{n})^T, T_n^{(i)} \rangle_n^2\langle X_{1}^{n} (X_{3}^{n})^T, T_n^{(i)} \rangle_n^2 \\
	& &&
	+4(n-1)\langle X_{1}^{n} (X_{1}^{n})^T, T_n^{(i)} \rangle_n\langle X_{1}^{n} (X_{2}^{n})^T, T_n^{(i)} \rangle_n^3 \Bigg] \\
	& \leq &&96\frac{(3n+1)nN_n}{(n-1)^4}K_1^2 K_2 \longrightarrow 0.
	\end{aligned}
	\end{equation}
	
	\paragraph{Term (I.4)}
	\begin{equation}
	\begin{aligned}
	& &&\frac{6N_n}{(n-1)^4} \sum_{i=1}^{N_n}\E\left[\left(\sum_{k=1}^n\langle \bar{X}^{n} (\bar{X}^{n})^T, T_n^{(i)} \rangle_n^2\right)^2\right] \\
	&= &&\frac{6N_n n^2}{(n-1)^4}\sum_{i=1}^{N_n}\E\left[\left(\langle \bar{X}^{n} (\bar{X}^{n})^T, T_n^{(i)} \rangle_n^2\right)^2\right] \\
	&\leq && \frac{6N_n n^2}{(n-1)^4} 8\Bigg(\underbrace{\sum_{i=1}^{N_n}\E\left[\left(\langle \bar{X}^{n} (\bar{X}^{n})^T - \frac{1}{n}\Sigma_n, T_n^{(i)} \rangle_n^2\right)^2\right]}_{=\frac{(n-1)^4}{n^4} \times \text{(term 5, Lemma 1)}} + \sum_{i=1}^{N_n}\E\left[\left(\langle \frac{1}{n}\Sigma_n, T_n^{(i)} \rangle_n^2\right)^2\right] \Bigg) \\
	& \leq &&\frac{6N_n n^2}{(n-1)^4} \times 8 \left( \frac{(n-1)^4}{n^4} \frac{M''_5}{n} K_1 K_2 + \frac{K_2}{n^4}\right) \longrightarrow 0,
	\end{aligned}
	\end{equation}
	using $M''_5$ from Equation \ref{eq_m5}.
	
	\paragraph{Term (I.5)}
	\begin{equation}
	\begin{aligned}
	& &&\frac{6N_n}{(n-1)^4} \sum_{i=1}^{N_n}\E\left[ \left(\sum_{k=1}^n 4\langle X_{k}^{n} (\bar{X}^{n})^T, T_n^{(i)} \rangle_n \langle \bar{X}^{n} (\bar{X}^{n})^T, T_n^{(i)} \rangle_n\right)^2\right] \\
	& = &&\frac{6\times 16N_n}{(n-1)^4} \sum_{i=1}^{N_n}\E\left[\left(\sum_{k=1}^n\langle \bar{X}^{n} (\bar{X}^{n})^T, T_n^{(i)} \rangle_n^2\right)^2\right] \longrightarrow 0.
	\end{aligned}
	\end{equation}
	
	\paragraph{Term (I.6)}
	\begin{equation}
	\begin{aligned}
	& &&\frac{6N_n}{(n-1)^4} \sum_{i=1}^{N_n}\E\left[\left(\sum_{k=1}^n2 \langle X_{k}^{n} (X_{k}^{n})^T, T_n^{(i)} \rangle_n \langle \bar{X}^{n} (\bar{X}^{n})^T, T_n^{(i)} \rangle_n\right)^2\right] \\
	& \leq &&\frac{12N_n}{(n-1)^4} \sum_{i=1}^{N_n} \left(\E\left[\left(\sum_{k=1}^n \langle X_{k}^{n} (X_{k}^{n})^T, T_n^{(i)} \rangle_n^2 \right)^2\right]  + \E\left[\left(\sum_{k=1}^n\langle \bar{X}^{n} (\bar{X}^{n})^T, T_n^{(i)} \rangle_n^2\right)^2\right] \right) \\
	& &&\longrightarrow 0.
	\end{aligned}
	\end{equation}
	
	The convergence to $0$ of each one of the 6 terms leads to the convergence of \textbf{term (I)}:
	\begin{equation}
	\begin{aligned}
	\V\left[\frac{1}{(n-1)^2}\sum_{i=1}^{N_n} \sum_{k=1}^n \langle \tilde{X}_{k}^{n} (\tilde{X}_{k}^{n})^T, T_n^{(i)} \rangle_n^2 \right] \longrightarrow 0.
	\end{aligned}
	\end{equation}
	
	\paragraph{Term (II)} We then show that $\V\left[\frac{1}{np_n}\sum_{i=1}^{N_n} \langle S_n T_n^{(i)}, T_n^{(i)} S_n \rangle_n \right] \longrightarrow 0$.
	\begin{equation}
	\begin{aligned}
		&\V\left[\frac{1}{np_n}\sum_{i=1}^{N_n} \langle S_n T_n^{(i)}, T_n^{(i)} S_n \rangle_n \right] 
		&&\leq \E\left[\frac{1}{n^2p_n^2} \left(\sum_{i=1}^{N_n} \langle S_n T_n^{(i)}, T_n^{(i)} S_n \rangle_n\right)^2 \right] \\
		& &&\leq \frac{N_n}{n^2p_n^2}\sum_{i=1}^{N_n}\E\left[\sum_{i=1}^{N_n} \langle S_n T_n^{(i)}, T_n^{(i)} S_n \rangle_n^2 \right]\\
		& &&\leq \frac{N_n^2}{n^2p_n} \E\left[\lVert S_n \rVert_n^2\right] \\
		& \V\left[\frac{1}{np_n}\sum_{i=1}^{N_n} \langle S_n T_n^{(i)}, T_n^{(i)} S_n \rangle_n \right] 
		&&\leq \frac{N_n^2}{n^2p_n} \E\left[\lVert S_n - \Sigma_n \rVert_n^2 + \lVert \Sigma_n \rVert_n^2\right].
	\end{aligned}
	\end{equation}
	Thanks to Lemma 1 \cite{Oriol2023}, particularly in the proof of the Lemma 1 equations (16) and (17) about $\beta_n^2 = \E\left[\lVert S_n - \Sigma_n \rVert_n^2\right]$, we have that $\frac{n}{p_n}\E\left[\lVert S_n - \Sigma_n \rVert_n^2\right]$ is bounded. Additionnally, $\lVert S_n \rVert_n^2 \leq \frac{p_n}{n}\sqrt{K_n}$. So:
	\begin{equation}
	\begin{aligned}
		&\V\left[\frac{1}{np_n}\sum_{i=1}^{N_n} \langle S_n T_n^{(i)}, T_n^{(i)} S_n \rangle_n \right] \longrightarrow 0.
	\end{aligned}
	\end{equation}
	
	\paragraph{Term (III)} We finally show that $\V\left[\frac{1}{n}\sum_{i=1}^{N_n} \langle S_n, T_n^{(i)} \rangle_n^2 \right] \longrightarrow 0$.
	\begin{equation}
	\begin{aligned}
		&\V\left[\frac{1}{n}\sum_{i=1}^{N_n} \langle S_n, T_n^{(i)} \rangle_n^2 \right] \leq \frac{N_n}{n^2} \sum_{i=1}^{N_n} \E\left[ \langle S_n, T_n^{(i)} \rangle_n^4 \right] \\
		&\V\left[\frac{1}{n}\sum_{i=1}^{N_n} \langle S_n, T_n^{(i)} \rangle_n^2 \right] \leq \frac{N_n}{n^2} \E\left[\lVert S_n \rVert_n^4 \right] \\
	\end{aligned}
	\end{equation}
	Thanks to the proof of Lemma 3 \cite{Oriol2023}, equation (54), we have that $\V\left[\lVert S_n \rVert_n^2 \right] $ is bounded. Moreover, from Lemma 1 \cite{Oriol2023}, we have that $\E\left[\lVert S_n \rVert_n^2 \right] $ is bounded too. So $\E\left[\lVert S_n \rVert_n^4 \right] $ is bounded. Consequently,
	\begin{equation}
	\begin{aligned}
		&\V\left[\frac{1}{n}\sum_{i=1}^{N_n} \langle S_n, T_n^{(i)} \rangle_n^2 \right]  \longrightarrow 0.
	\end{aligned}
	\end{equation}
	
	Backing up, the convergences of \textbf{Term (I)-(III)} prove that:
	\begin{equation}
	\begin{aligned}
		&\E\left[\left(\sum_i \hat{\V}(\langle S_n, T_n^{(i)} \rangle_n) - \V(\langle S_n, T_n^{(i)} \rangle_n)\right)^2\right] \longrightarrow 0.
	\end{aligned}
	\end{equation}
	
	\subsection{Proof of Theorem 2}
	We denote:
	\begin{equation}
	\begin{aligned}
	&a_n^2 = c_1 \det(A_n), \alpha_n = c_1^* \det(A_n), d_n^2 = \det(A_n), \delta_n^2 = \E\left[\left\lVert S_n - \sum_{k} \langle \Sigma_n, T_n^{(k)} \rangle_n T_n^{(k)} \right\rVert_n^2\right], \\
	&u_n^2 = (a_n^2 - \alpha_n)^2.
	\end{aligned}
	\end{equation}
	We will first prove that $\E\left[\left \lVert \tilde{S}_n^* - \tilde{\Sigma}_n^{*}\right \rVert_n^2\right] \rightarrow 0$ and deduce that $\E\left[\left \lVert S_n^* - \Sigma_n^{*}\right \rVert_n^2\right] \rightarrow 0$ by a property of the orthogonal projection.  \\
	We have:
	\begin{equation}
	\begin{aligned}
	\E\left[\left \lVert \tilde{S}_n^* - \tilde{\Sigma}_n^{*}\right \rVert_n^2\right] &= \E\left[\frac{u_n^2}{d_n^2} \right] + \sum_k \V\left[\langle S_n, T_n^{(k)} \rangle_n\right].
	\end{aligned}
	\end{equation}
	Firstly, we have thanks to Lemma 1: 
	\begin{equation}\label{sqrt}
	\begin{aligned}
	\sum_k \V\left[\langle S_n, T_n^{(k)} \rangle_n\right] &= \sum_k \E\left[\langle S_n - \Sigma_n, T_n^{(k)} \rangle_n^2 \right] \\
	&\leq \sum_k \sqrt{\E\left[\langle S_n - \Sigma_n, T_n^{(k)} \rangle_n^4 \right]} \\
	&\leq \sqrt{N_n \sum_k \E\left[\langle S_n - \Sigma_n, T_n^{(k)} \rangle_n^4 \right]} \\
	&\leq \sqrt{\frac{N_nMK_1K_2}{n}}\\
	\sum_k \V\left[\langle S_n, T_n^{(k)} \rangle_n\right] &\longrightarrow 0.
	\end{aligned}
	\end{equation}
	Then, we have that $d_n^2 - \delta_n^2 \underset{q.m}{\longrightarrow} 0$. Indeed:
	\begin{equation}
	\begin{aligned}
	d_n^2 - \delta_n^2 = \lVert S \rVert_n^2 - \E[\lVert S \rVert_n^2] - \sum_k \left(\langle \Sigma_n, T_n^{(k)} \rangle_n^2-\langle S_n, T_n^{(k)} \rangle_n^2 \right)
	\end{aligned}
	\end{equation}
	So, we have:
	\begin{equation}
	\begin{aligned}
	\E\left[\left(d_n^2 - \delta_n^2\right)^2\right] &\leq 2\V\left[\lVert S \rVert_n^2\right] + 2 \V\left[\sum_k \langle \Sigma_n-  S_n, T_n^{(k)} \rangle_n^2\right]\\
	&\leq 2\V\left[\lVert S \rVert_n^2\right] + 2 \E\left[\left(\sum_k \langle \Sigma_n-  S_n, T_n^{(k)} \rangle_n^2\right)^2\right]\\
	&\leq 2\V\left[\lVert S \rVert_n^2\right] + 2 N_n\sum_k \E\left[\langle \Sigma_n-  S_n, T_n^{(k)} \rangle_n^4\right]\\
	\E\left[\left(d_n^2 - \delta_n^2\right)^2\right]&\leq 2\V\left[\lVert S \rVert_n^2\right] + 2 N_n\frac{MK_1K_2}{n}. &&\text{(Lemma 1)}\\
	\end{aligned}
	\end{equation}
	$\V\left[\lVert S \rVert_n^2\right] \rightarrow 0$ from the proof of Lemma 3 in \cite{Oriol2023} and from the proof of Lemma 3.3 in \cite{Ledoit2004}. So we have $d_n^2 - \delta_n^2 \underset{q.m}{\longrightarrow} 0$.  \\
	Now, let's prove $a_n^2 - \E[\alpha_n]  \underset{q.m}{\longrightarrow} 0$.  \\
	Firstly, let's remark that $0 \leq \E[\alpha_n] \leq \delta_n^2$. Indeed:
	\begin{equation}
	\begin{aligned}
	\E[\alpha_n] &= \E\left[\langle S_n, \Sigma_n \rangle_n - \sum_k \langle S_n, T_n^{(k)} \rangle_n\langle \Sigma_n, T_n^{(k)} \rangle_n \right]  \\
	&= \lVert \Sigma_n \rVert_n^2 - \sum_k \langle \Sigma_n, T_n^{(k)} \rangle_n^2 \\
	&= \left\lVert \Sigma_n- \sum_k \langle \Sigma_n, T_n^{(k)} \rangle_nT_n^{(k)} \right \rVert_n^2 \\
	\E[\alpha_n]& \geq 0
	\end{aligned}
	\end{equation}
	And,
	\begin{equation}
	\begin{aligned}
	\delta_n^2 &= \E\left[\left\lVert S_n - \sum_{k} \langle \Sigma_n, T_n^{(k)} \rangle_n T_n^{(k)} \right\rVert_n^2\right]\\
	&= \E\left[\lVert S_n \rVert_n^2\right] - \sum_{k} \langle \Sigma_n, T_n^{(k)} \rangle_n^2  \\
	&= \E[\alpha_n] + \E\left[\lVert S_n - \Sigma_n \rVert_n^2\right]\\
	\delta_n^2& \geq \E[\alpha_n] 
	\end{aligned}
	\end{equation}
	So, we have $0 \leq \E[\alpha_n] \leq \delta_n^2$.  \\
	Let's find a lower and an upper bound for $a_n^2 - \E[\alpha_n]$. For the upper bound:
	\begin{equation}
	\begin{aligned}
	a_n^2 - \E[\alpha_n] &= \min((\hat{c_1}d_n^2)_+,d_n^2) -  \E[\alpha_n]  \\
	& \leq (\hat{c_1}d_n^2)_+ -   \E[\alpha_n]\\
	& \leq | \hat{c_1}d_n^2 -  \E[\alpha_n]| && \text{(as }\E[\alpha_n] \geq 0 \text{)}  \\
	a_n^2 - \E[\alpha_n]& \leq \max\left(|\hat{c_1}d_n^2 -  \E[\alpha_n]|, |d_n^2 - \delta_n^2| \right)
	\end{aligned}
	\end{equation}
	For the lower bound:
	\begin{equation}
	\begin{aligned}
	a_n^2 -  \E[\alpha_n]&= \min((\hat{c_1}d_n^2)_+-  \E[\alpha_n],d_n^2-  \E[\alpha_n])   \\
	&\geq \min((\hat{c_1}d_n^2)_+-  \E[\alpha_n],d_n^2-  \delta_n^2) && (\text{as } \E[\alpha_n] \leq \delta_n^2)   \\
	&\geq \min(-|\hat{c_1}d_n^2-  \E[\alpha_n]|,-| d_n^2-  \delta_n^2|)   \\
	a_n^2 -  \E[\alpha_n]& \geq -\max\left(|\hat{c_1}d_n^2 - \alpha_n|, |d_n^2 - \delta_n^2| \right)
	\end{aligned}
	\end{equation}
	So,
	\begin{equation}
	\begin{aligned}
	\E\left[\left(a_n^2 -  \E[\alpha_n]\right)^2\right]&\leq \E\left[\max\left(|\hat{c_1}d_n^2 - \alpha_n|, |d_n^2 - \delta_n^2| \right)^2 \right]  \\
	& \leq \E\left[\left(\hat{c_1}d_n^2 - \alpha_n\right)^2\right] + \E\left[\left(d_n^2 - \delta_n^2 \right)^2 \right]
	\end{aligned}
	\end{equation}
	From Corollary 1, $\E\left[\left(\hat{c_1}d_n^2 - \alpha_n\right)^2\right] \rightarrow 0$, and we proved above that $\E\left[\left(d_n^2 - \delta_n^2 \right)^2 \right] \rightarrow 0$. So, $a_n^2 - \E[\alpha_n]  \underset{q.m}{\longrightarrow} 0$.
	Then, let's prove that $\V\left[\alpha_n\right] \rightarrow 0$.
	\begin{equation}
	\begin{aligned}
	\V[\alpha_n] &= \V[c_1^* \det(A_n)]\\
	&= \V\left[\langle S_n, \Sigma_n \rangle_n - \sum_k \langle S_n, T_n^{(k)} \rangle_n\langle \Sigma_n, T_n^{(k)} \rangle_n \right]\\
	\V[\alpha_n]& \leq 2\V\left[\langle S_n, \Sigma_n \rangle_n\right] + 2\V\left[\sum_k \langle S_n, T_n^{(k)} \rangle_n\langle \Sigma_n, T_n^{(k)} \rangle_n \right].\\
	\end{aligned}
	\end{equation}
	$\V\left[\langle S_n, \Sigma_n \rangle_n\right] \rightarrow 0$ from the proof of Theorem 3 \cite{Oriol2023} for the unknown mean and the proof of Theorem 3.3 \cite{Ledoit2004} for the known mean.  \\
	For the second term, we have:
	\begin{equation}
	\begin{aligned}
	&\V\left[\sum_k \langle S_n, T_n^{(k)} \rangle_n\langle \Sigma_n, T_n^{(k)} \rangle_n \right] \\
	&= \E\left[\left(\sum_k \langle S_n, T_n^{(k)} \rangle_n\langle \Sigma_n, T_n^{(k)} \rangle_n\right)^2 \right] - \left(\sum_k \langle \Sigma_n, T_n^{(k)} \rangle_n^2\right)^2\\
	& \leq \E\left[\sum_k \langle S_n, T_n^{(k)} \rangle_n^2\right]\sum_k \langle \Sigma_n, T_n^{(k)} \rangle_n^2  - \left(\sum_k \langle \Sigma_n, T_n^{(k)} \rangle_n^2\right)^2\\
	& = \left(\sum_k\V\left[ \langle S_n, T_n^{(k)} \rangle_n\right]\right)\left(\sum_k \langle \Sigma_n, T_n^{(k)} \rangle_n^2\right)\\
	& \leq \sqrt{\frac{N_nMK_1K_2}{n}} \left(\sum_k \langle \Sigma_n, T_n^{(k)} \rangle_n^2\right) \text{  (Eq \ref{sqrt})}\\
	& \leq \sqrt{\frac{N_nMK_1K_2}{n}} \lVert \Sigma_n \rVert^2\\
	& \longrightarrow 0.\\
	\end{aligned}
	\end{equation}
	So, as both terms converges to $0$, we have that $\V[\alpha_n] \rightarrow 0$.  \\
	Easily, we have then that $\E[u_n^2] \longrightarrow 0$ as $\E[u_n^2] \leq 2\left(\E[(a_n^2 - \E[\alpha_n])^2] + \V[\alpha_n] \right)$. \\
	Moreover, let's check that $0 \leq \frac{u_n^2}{d_n^2} \leq 2(d_n^2 + \delta_n^2)$.
	\begin{equation}
	\begin{aligned}
	|\alpha_n| &= \left|\langle S_n, \Sigma_n \rangle_n - \sum_k \langle S_n, T_n^{(k)} \rangle_n\langle \Sigma_n, T_n^{(k)} \rangle_n\right|\\
	&= \left|\left \langle S_n - \sum_k \langle S_n, T_n^{(k)} \rangle_nT_n^{(k)}, \Sigma_n - \sum_k \langle \Sigma_n, T_n^{(k)} \rangle_nT_n^{(k)}\right \rangle_n\right| \\
	& \leq \left \lVert S_n - \sum_k \langle S_n, T_n^{(k)} \rangle_nT_n^{(k)} \right\rVert_n \left \lVert \Sigma_n - \sum_k \langle \Sigma_n, T_n^{(k)} \rangle_nT_n^{(k)}\right \rVert_n\\
	|\alpha_n| & \leq d_n \delta_n.
	\end{aligned}
	\end{equation}
	And, by construction, $0 \leq a_n^2 \leq d_n^2$, so $0 \leq \frac{u_n^2}{d_n^2} \leq 2(d_n^2 + \delta_n^2)$.  \\
	Backing up, we proved that $d_n^2 - \delta_n^2 \underset{q.m}{\longrightarrow} 0$, $\E[u_n^2] \longrightarrow 0$, and $0 \leq \frac{u_n^2}{d_n^2} \leq 2(d_n^2 + \delta_n^2)$.
	We then can apply Lemma A.1 from \cite{Ledoit2004} with $\tau_1=2, \tau_2 = 0$, which proves that $\E\left[\frac{u_n^2}{d_n^2} \right] \longrightarrow 0$. \\
	Backing up again, it proves that $\tilde{S}_n^* - \tilde{\Sigma}_n^{*} \underset{q.m}{\longrightarrow} 0$.  \\
	Finally, $\mathcal{P}_{\mathcal{S}_n^+}$ is 1-Lipschitz for the norm $\lVert \cdot \rVert_n$, so $S_n^* - \Sigma_n^{*} \underset{q.m}{\longrightarrow} 0$.  \\
	The second part of the theorem follows:
	\begin{equation}
	\begin{aligned}
	\E\left[\left|\left\lVert S_n^* - \Sigma_n \right \rVert_n^2 - \left\lVert {\Sigma}_n^{*} - \Sigma_n \right \rVert_n^2 \right|\right] & = \E\left[\left|\left\langle {S}_n^* - {\Sigma}_n^{*}, {S}_n^* +{\Sigma}_n^{*} - 2\Sigma_n\right\rangle_n\right|\right]  \\
	& \leq \sqrt{\E\left[\lVert {S}_n^* - {\Sigma}_n^{*}\rVert_n^2\right]} \sqrt{\E\left[\lVert {S}_n^* + {\Sigma}_n^{*} - 2\Sigma_n \rVert_n^2\right]}.
	\end{aligned}
	\end{equation}
	Moreover, $\E\left[\lVert {S}_n^* - \Sigma_n \rVert_n^2\right]$ is bounded. In fact:
	\begin{equation}
	\begin{aligned}
	\E\left[\lVert {S}_n^* - \Sigma_n \rVert_n^2\right] & \leq \E\left[\lVert \tilde{S}_n^* - \Sigma_n \rVert_n^2\right] \\
	&\leq \E\left[2\lVert S_n - \Sigma_n \rVert_n^2 + 2 \left\lVert \sum_k \langle S_n, T_n^{(k)} \rangle_n  T_n^{(k)} - \Sigma_n \right\rVert_n^2  \right]  \\
	\E\left[\lVert {S}_n^* - \Sigma_n \rVert_n^2\right]& \leq 2\E\left[\lVert S_n - \Sigma_n \rVert_n^2 \right] + 2\lVert \Sigma_n \rVert_n^2 + 2\sum_k \V\left[ \langle S_n, T_n^{(k)} \rangle_n  \right]
	\end{aligned}
	\end{equation}
	$\E\left[\lVert S_n - \Sigma_n \rVert_n^2 \right]$ and $\lVert \Sigma_n \rVert_n^2$ are bounded from Lemma 3.1 in \cite{Ledoit2004} and Lemma 1 in \cite{Oriol2023}, and $\sum_k \V\left[ \langle S_n, T_n^{(k)} \rangle_n  \right] \rightarrow 0$, so $\E\left[\lVert {S}_n^* - \Sigma_n \rVert_n^2\right]$ is bounded.
	So, using the first part of the Theorem 2 we've just proved, we have that: $\E\left[\left|\left\lVert {S}_n^* - \Sigma_n \right \rVert_n^2 - \left\lVert {\Sigma}_n^{*} - \Sigma_n \right \rVert_n^2 \right|\right] \rightarrow 0$.

\end{appendices}

\bibliography{mybib}

\end{document}